# PINNING OF POLYMERS AND INTERFACES BY RANDOM POTENTIALS

By Kenneth S. Alexander[1] and Vladas Sidoravicius

*University of Southern California and IMPA*

We consider a polymer, with monomer locations modeled by the trajectory of a Markov chain, in the presence of a potential that interacts with the polymer when it visits a particular site 0. Disorder is introduced by, for example, having the interaction vary from one monomer to another, as a constant $u$ plus i.i.d. mean-0 randomness. There is a critical value of $u$ above which the polymer is pinned, placing a positive fraction of its monomers at 0 with high probability. This critical point may differ for the quenched, annealed and deterministic cases. We show that self-averaging occurs, meaning that the quenched free energy and critical point are nonrandom, off a null set. We evaluate the critical point for a deterministic interaction ($u$ without added randomness) and establish our main result that the critical point in the quenched case is strictly smaller. We show that, for every fixed $u \in \mathbb{R}$, pinning occurs at sufficiently low temperatures. If the excursion length distribution has polynomial tails and the interaction does not have a finite exponential moment, then pinning occurs for all $u \in \mathbb{R}$ at arbitrary temperature. Our results apply to other mathematically similar situations as well, such as a directed polymer that interacts with a random potential located in a one-dimensional defect, or an interface in two dimensions interacting with a random potential along a wall.

**1. Introduction.** Consider the following physical problems.

Problem 1.1. A polymer molecule in $d$ dimensions is in the presence of a potential well, at a single site, which attracts the monomers located there. The configuration of the polymer—that is, the sequence of locations of the monomers—follows the trajectory of a random walk. The polymer can configure itself so that a positive fraction of all monomers are at the

Received July 2005.

[1]Supported in part by NSF Grants DMS-01-03790 and DMS-04-05915.

*AMS 2000 subject classifications.* Primary 82D60; secondary 82B44, 60K35.

*Key words and phrases.* Pinning, polymer, disorder, interface, random potential.

potential well, in which case we say the polymer is *pinned*, but of course this involves a sacrifice of entropy. For what values of potential depth and temperature will the polymer be pinned? If there is a depinning transition, what are the critical exponents that describe it? We can introduce disorder by allowing the polymer to be a *heteropolymer*, meaning that the attraction to the potential well varies from one monomer to another, perhaps being a repulsion for some monomers. This problem is examined in [15].

Problem 1.2. Consider the two-dimensional Ising model below the critical temperature in a square box with plus boundary condition on the bottom side (the wall) and minus boundary condition on the other three sides. This forces the existence of an interface connecting the lower left and lower right corners of the box. Suppose that the interaction between wall sites and adjacent box sites is weaker than the bulk Ising interaction, say, $(1-u)J$ instead of $J$; this gives the interface an energetic advantage for each location where it touches the wall. For what values of $u$ and temperature will the interface be pinned to the wall? The absence of pinning is called *wetting*; when there is a transition from pinning to wetting, what critical exponents describe it? Again we can introduce disorder by allowing the value of $u$ to vary from site to site along the wall. This problem is examined in [1, 5, 6, 8, 9]. To fit it in our present context, we must impose an solid-on-solid (SOS) restriction, as discussed below.

Problem 1.3. Consider a polymer in $d + 1$ dimensions, directed in one coordinate, in the presence of a wall that confines the polymer to a halfspace, and suppose that a potential attracts those monomers that touch the wall. As with the potential well, we may ask, when is the polymer pinned, and what is the nature of the depinning transition, if any? Such questions arise in the study of adhesion. Disorder may be introduced by considering heteropolymers or by allowing the potential to vary randomly from site to site on the wall. This problem is considered in [19, 20, 24].

Problem 1.4. Consider a polymer in $d + 1$ dimensions, directed in the $(d + 1)$st coordinate, in the presence of a lower-dimensional *defect*, meaning a subspace where a potential attracts those monomers located in it. We take the defect to be the coordinate axis for the $(d + 1)$st coordinate. In superconducting materials, under certain conditions, nearly all magnetic flux lines are confined to a small number of random tubes, and the trajectories of these tubes are influenced by defects which may attract them. The "polymer" here has been used to model such trajectories. Superconducting properties are affected by whether the defects pin the flux tubes. Again, disorder may be introduced by way of either a heteropolymer or site-to-site variation in the defect; these are mathematically equivalent. In fact, if we view the $(d + 1)$st coordinate as merely an index for the monomers of a polymer existing in $d$ dimensions, we see that Problem 1.1 is equivalent as well. This problem is considered in [10, 11, 17].

All of these problems (with the exception of the randomly varying potential in a wall in Problem 1.3, when $d \geq 2$) share the following mathematical setup. There is an underlying Markov chain $\{X_i, i \geq 0\}$ representing the trajectories in the absence of the potential; it is governed by a transition probability $p(\cdot, \cdot)$ and has a state space $\Sigma$. For convenience, we label the directed coordinate as "time" so that the polymer trajectories become space–time trajectories of the chain, and assume the chain is irreducible and aperiodic. There is a unique site in $\Sigma$ which we call 0 where the potential is located; we consider trajectories of length $n$ starting from state 0 at time 0 and denote the corresponding measure $P^X_{[0,n]}$. The potential at 0 at time $i$ has form $u + V_i$, where the $V_i$ are i.i.d. with mean zero; we refer to $\{V_i : i \geq 1\}$ as the *disorder*. For $n$ and a realization $\{V_i, i \geq 1\}$ fixed, we attach a Gibbs weight

$$(1.1) \qquad \exp\left(\beta \sum_{i=1}^n (u + V_i) \delta_{\{x_i = 0\}}\right) P^X_{[0,n]}(\{x_i, 0 \leq i \leq n\})$$

to each trajectory $\{x_i : 0 \leq i \leq n\}$. Here $\delta_A$ denotes the indicator of an event $A$. The corresponding partition function and finite-volume Gibbs distribution on length-$n$ trajectories are denoted $Z^{\{V_i\}}_{[0,n]}(\beta, u)$ and $\mu^{\beta, u, \{V_i\}}_{[0,n]}$, respectively. We omit the $\{V_i\}$ when $V_i \equiv 0$. We write $P^V_{[a,b]}$ for the distribution of $(V_a, V_{a+1}, \ldots, V_b)$ and $P^V$ for the distribution of the full sequence $\{V_i\}$. We use $\langle \cdot \rangle$ to denote expected value, with superscripts and subscripts corresponding to the measure, so that, for example, $\langle \cdot \rangle^V$ denotes expectation under $P^V$. Let

$$L_n = \sum_{i=1}^n \delta_{\{X_i = 0\}}.$$

For fixed $\beta, u$, we say the polymer is *pinned* at $(\beta, u)$ if, for some $\delta > 0$,

$$\lim_n \mu^{\beta, u, \{V_i\}}_{[0,n]}(L_n > \delta n) = 1, \qquad P^V_{[0,\infty)}\text{-a.s.}$$

It is clear that if the polymer is pinned at $(\beta, u)$, then it is pinned at $(\beta, u')$ for all $u' > u$. Therefore, there is a (possibly infinite) critical $u_c(\beta, \{V_i\})$ such that the polymer is pinned for $u > u_c(\beta, \{V_i\})$ and not pinned for $u < u_c(\beta, \{V_i\})$. In Theorem 3.1 we will establish that self-averaging holds for the free energy, which implies that there is a nonrandom *quenched critical point* $u_c^q = u_c^q(\beta)$ such that $u_c(\beta, \{V_i\}) = u_c^q(\beta)$ with $P^V$-probability one.

There are two other critical points to consider. The *deterministic critical point* $u_c^d = u_c^d(\beta)$ is the critical point for the *deterministic model*, which is the case $V_i \equiv 0$. The *annealed model* is obtained by averaging the Gibbs weight (1.1) over the disorder; the annealed model at $(\beta, u)$ is thus the same as the deterministic model at $(\beta, u + \beta^{-1} \log M_V(\beta))$, where $M_V$ is the moment generating function of $V_1$, and the corresponding *annealed critical point* is $u_c^a = u_c^a(\beta) = u_c^d(\beta) - \beta^{-1} \log M_V(\beta)$. It is not hard to show that

$$(1.2) \qquad\qquad u_c^a \leq u_c^q \leq u_c^d;$$

in fact, once we establish the existence of the quenched free energy (Theorem 3.1), the first inequality is an immediate consequence of Jensen's inequality. It is the strictness of the second inequality that is less obvious, and that we establish here (Theorem 1.5).

The relation between these critical points may be interpreted heuristically as follows. Let $E_i$ denote the (possibly infinite) $i$th excursion length for the chain, that is, the time between the $(i-1)$st and $i$th visits to 0, with the visit at time 0 counted as the 0th visit. We consider for this heuristic the case in which $P^X(n \leq E_1 < \infty)$ does not decay exponentially, as the heuristics are somewhat different for the alternative. For $M \in \mathbb{R}$ and $\delta \in [0,1]$, consider trajectory/disorder pairs $(\{x_i\}, \{v_i\})$ for which the fraction of time at 0 is approximately $\delta$, and the average random potential experienced there is approximately $M$:

$$L_n \approx \delta n, \qquad \frac{1}{L_n} \sum_{i=1}^{n} v_i \delta_{\{x_i=0\}} \approx M.$$

The annealed partition function

$$Z_{[0,n]}^a(\beta, u) = \sum_{\text{walks } \{x_i\}} \left\langle \exp\left( \beta \sum_{i=1}^{n} (u + V_i) \delta_{\{x_i=0\}} \right) \right\rangle_{[1,n]}^V P_{[0,n]}^X(\{x_i, 0 \leq i \leq n\})$$

can be decomposed into contributions from various $M$ and $\delta$. From basic large deviation theory, the log of such a contribution is asymptotically

$$\left( \beta(u + M) - I_V(M) - I_E\left(\frac{1}{\delta}\right) \right) \delta n,$$

where $I_V$ and $I_E$ are large-deviation rate functions for $V_1$ and $E_1$, respectively:

$$I_V(M) = -\lim_{\varepsilon \searrow 0} \lim_n \frac{1}{n} \log P_{[1,n]}^V\left( \frac{1}{n} \sum_{i=1}^{n} V_i \in (M - \varepsilon, M + \varepsilon) \right)$$

$$= \sup_x (tx - \log\langle e^{xV_1}\rangle^V)$$

and

$$I_E(t) = -\lim_{\varepsilon \searrow 0} \lim_n \frac{1}{n} \log P^X\left( \frac{1}{n} \sum_{i=1}^{n} E_i \in (t - \varepsilon, t + \varepsilon) \right),$$

where the supremum is over all $x$ with $\langle e^{xV_1} \rangle^V < \infty$. The annealed free energy $f^a(\beta, u)$ is then given by

$$
\begin{aligned}
(1.3) \qquad \beta f^a(\beta, u) &= \lim_n \frac{1}{n} \log Z^a_{[0,n]}(\beta, u) \\
&= \sup_{M, \delta} \Big( \beta(u + M) - I_V(M) - I_E\Big(\frac{1}{\delta}\Big) \Big) \delta.
\end{aligned}
$$

Since the free energy associated to "unpinned" trajectories is 0 (see Theorem 2.1 for a precise statement), we expect pinning to occur precisely when the free energy is positive, that is, when the contribution from some $M, \delta$ outweighs the contribution from unpinned trajectories:

$$
\beta(u + M) - I_V(M) - I_E\Big(\frac{1}{\delta}\Big) > 0 \qquad \text{for some } M \in \mathbb{R}, \delta > 0.
$$

(The jump from a free energy statement to a pathwise statement, i.e., pinning, is not trivial here; see [3] for a rigorous derivation in a related context. Our derivation in this introduction is heuristic only.) It is easy to see that, since we are assuming $E_1$ has no finite exponential moment, we have $I_E(1/\delta) \searrow -\log P^X(E_1 < \infty)$ as $\delta \to 0$, while, from basic large deviation theory, we have

$$
\sup_M (\beta M - I_V(M)) = \log M_V(\beta).
$$

Therefore,

$$
u^a_c(\beta) = -\beta^{-1}(\log M_V(\beta) + \log P^X(E_1 < \infty)),
$$

which in the deterministic case says that

$$
u^d_c(\beta) = -\beta^{-1} \log P^X(E_1 < \infty).
$$

The reason the heuristic does not apply in the quenched case, meaning we need not have $u^q_c = u^a_c$, is that, for a fixed realization $\{v_i\}$ of $\{V_i\}$, the sample

$$
\{v_i : X_i = 0\}
$$

of potentials selected out of $\{v_i : 1 \le i \le n\}$ by the chain via its return times to 0 is not an i.i.d. sample from $P^V$, so the large deviation rate $I_V$ does not apply. For one thing, in the annealed case, when the chain selects $\delta n$ potentials $V_i$ averaging to some $M > 0$, it means that, with high probability, the overall disorder of $n$ potentials is very atypical, averaging about $M\delta$ though its mean is 0. If instead we have a "typical" disorder, averaging to near 0, the cost is greater (i.e., the probability is lower) to select a size $\delta n$ sample with average $M$. Further, the chain selects the sample from the realization without replacement, which again increases the cost of large deviations. Overall, compared to an i.i.d. sample, the chain achieves large deviation averages at greater, or at best equal, cost by its returns-to-0 sampling procedure, for a "typical" fixed realization.

If we find a case in which $u_c^q = u_c^a$, then, we may interpret this as meaning that the chain is "almost as efficient as an i.i.d. selection" in obtaining large-deviation averages via its sampling procedure. Similarly, the weaker statement that $u_c^q < u_c^d$ [i.e., strict inequality in (1.2)] means that the chain can obtain large-deviation averages at low enough cost that when the mean of $u + V_i$ is slightly too small to induce pinning, that is, $u < u_c^d$, the chain can compensate, without excessive cost, by returning to 0 at on-average-favorable times. Our main result here is that this weaker statement is true in great generality. In the physics literature, the belief, based mainly on nonrigorous methods [5], analogous to periodic potentials [16] and numerics [15], is that the stronger statement $u_c^q = u_c^a$ is not always true.

In Problems 1.1 and 1.4 the underlying chain is a symmetric simple random walk on $\mathbb{Z}^d$. For $d = 1, 2$, the deterministic critical point is well known to be 0, so the conclusion $u_c^q < 0$ means that, even when the disorder is on average slightly negative, the chain will be pinned. This result was obtained in [10] for a periodic potential, which is frequently used in the physics literature as a surrogate for a random one.

Fitting Problem 1.2, on the Ising interface, into our setup requires some tweaking. First, we must impose the standard solid-on-solid (SOS) restriction, meaning overhangs in the interface are prohibited. In a box $[-L, L] \times [0, 2L]$, the interface is then described by a sequence of nonnegative integer heights $x_{-L}, \ldots, x_L$, with $x_i = m$ meaning that the interface above site $i$ is between $m - 1$ and $m$. Second, we must consider only the energetic cost of the interface itself, which is twice its length when the interaction is equal everywhere, and not consider the effect of the interface on the partition functions for the regions above and below it. We ignore horizontal bonds in calculating the length of an interface, since every allowed interface has the same number $2L + 1$ of them. Thus, the Gibbs weight of an interface is

$$\exp\left(2\beta\left(\sum_{i=-L}^{L}(u + V_i)\delta_{\{x_i=0\}} - \sum_{i=-L}^{L+1}|x_i - x_{i-1}|\right)\right).$$

Taking $u = 0, V_i \equiv 0$, we see that the underlying Markov chain is a random walk with transition probability $p(m, n)$ proportional to $\exp(-2\beta|n - m|)$, conditioned to stay nonnegative and to be 0 at times $-L - 1$ and $L + 1$.

The model (1.1), with symmetric simple random walk on $\mathbb{Z}^d$ as the underlying Markov chain, is also related to a special case of the Anderson model (on a lattice) in which the potential is nonzero at just a single site. In the nonrandom case $V_i \equiv 0$, pinning in (1.1) corresponds to localization in the Anderson model. For a random potential and discrete time, letting $u_A(n, x)$ denote the contribution to the partition function $Z_{[0,n]}^{\{V_i\}}(\beta, u)$ from paths ending at $x$ at time $n$, it is easily seen that $u_A$ satisfies

$$u_A(n+1, x) - u_A(n, x) = \tfrac{1}{2}e^{\beta V_n \delta_0(x)} \Delta u_A(n, x) + (e^{\beta V_n \delta_0(x)} - 1)u_A(n, x),$$

where $\Delta$ denotes the discrete Laplacian in the space variable $x$. This is a discrete time analog of the continuous-time equation seen in the corresponding case of the Anderson model:

$$du_A(t, x) = \Delta u_A(t, x)\, dt + \beta u_A(t, x)\delta_0(x)\, dW(t),$$

where $W$ is Brownian motion, meaning the disorder is a white noise in time, at 0. See [4] and [12] for more on the Anderson model.

For a trajectory $\{x_i, i \geq 0\} \in \Sigma^\infty$, an *excursion* (from 0) is either a segment $\{x_i, s \leq i \leq t\}$ with $x_s = x_t = 0, x_i \neq 0$ for $s < i < t$, or a segment $\{x_i, s \leq i < \infty\}$ with $x_s = 0, x_i \neq 0$ for $i > s$. The length of the excursion is $t - s$ or $\infty$ respectively, and we write $\mathcal{E}_i$ for the $i$th excursion and $E_i$ for its length, when $L_\infty \geq i$. When $L_\infty = i$, we further define $E_{i+1}$ to be $\infty$.

Our main result is the following. The proof is in Section 4.

THEOREM 1.5.  *Let $\{X_i\}$ be an irreducible Markov chain with state space containing a state 0, and let the potential $\{V_i\}$ at 0 be i.i.d. with mean 0, positive variance and $\langle e^{t|V_1|}\rangle^V < \infty$ for some $t > 0$. Suppose $u_c^d(1) > -\infty$. Then the pinning critical points for the measures $\mu_{[0,n]}^{\beta, u, \{V_i\}}$ and $\mu_{[0,n]}^{\beta, u}$ satisfy $-\infty < u_c^q(\beta) < u_c^d(\beta) < \infty$ for all $\beta > 0$. In particular, there is pinning at $(\beta, u_c^d(\beta))$.*

From Theorem 2.1 below, to satisfy the condition $u_c^d(1) > -\infty$ in Theorem 1.5, it is necessary and sufficient that either $E_1$ has no finite exponential moment, or the interval where the moment generating function of $E_1$ is finite includes its upper endpoint.

Define

$$G_V(x) = P^V(V_1 < x),$$

$$\overline{G}_V(x) = 1 - G_V(x),$$

$$\overline{G}_V^{-1}(t) = \sup\{x \geq 0 : \overline{G}_V(x) \geq t\},$$

with $\overline{G}_V^{-1}(t)$ defined to be 0 if $t > \overline{G}_V(0)$.

Theorem 1.5 says that the critical curve $u = u_c^q(\beta)$ in the $(\beta, u)$ plane is strictly below the curve $u = u_c^d(\beta)$. Further information about this curve is contained in the following theorem. The theme here is that when the disorder distribution has a sufficiently fat positive tail, pinning can occur by virtue of the chain returning to 0 at only those times when the disorder is exceptionally large. Related behavior is considered in the physics literature in [22].

THEOREM 1.6.   *Let $\{X_i\}$ be an irreducible Markov chain with state space containing a site 0, and let the potential $\{V_i\}$ at 0 be i.i.d. and nonconstant, with mean 0:*

(i) *If $V_1$ is unbounded, then, for each $u \in \mathbb{R}$, there is pinning at $(\beta, u)$ for all sufficiently large $\beta$.*

(ii) *If there exists a subsequence $k_j \to \infty$ satisfying*

$$(1.4) \qquad \min_{k_0 \le k \le k_j} P^X(E_1 = k) \ge e^{-o(\overline{G}_V^{-1}(1/k_j))} \qquad \text{as } j \to \infty,$$

*then there is pinning at $(\beta, u)$ for all $\beta > 0, u \in \mathbb{R}$.*

COROLLARY 1.7.   *Let $\{X_i\}$ be an irreducible Markov chain with state space containing a site 0, and let the potential $\{V_i\}$ at 0 be i.i.d. with mean 0:*

(i) *If $P^X(E_1 = k) \ge Ck^{-\gamma}$ for all sufficiently large $k$ for some $C > 0, \gamma > 1$, and $V_1$ does not have a finite exponential moment, then there is pinning at $(\beta, u)$ for all $\beta > 0, u \in \mathbb{R}$.*

(ii) *If $P^X(E_1 = k) \ge Ce^{-\alpha k}$ for all sufficiently large $k$ for some $C, \alpha > 0$, and the positive part $V_1^+$ of $V_1$ satisfies $E^V((V_1^+)^\theta) = \infty$ for some $0 < \theta < 1$, then there is pinning at $(\beta, u)$ for all $\beta > 0, u \in \mathbb{R}$.*

(iii) *More generally, if there exist a decreasing positive function $p$ and an increasing positive function $w$ on $[0, \infty)$ satisfying*

$$P^X(E_1 = k) \ge p(k) \qquad \text{for all sufficiently large } k, \qquad \sum_{k=1}^{\infty} \frac{1}{w(k)} < \infty$$

(1.5)

*and, letting $\zeta(x) = \log 1/x$,*

$$(1.6) \qquad E^V((\zeta \circ p \circ w)^{-1}(\varepsilon V_1^+)) = \infty \qquad \text{for all } \varepsilon > 0,$$

*then there is pinning at $(\beta, u)$ for all $\beta > 0, u \in \mathbb{R}$. Here $\circ$ denotes composition, and $(\zeta \circ p \circ w)^{-1}(x)$ is defined to be 0 if $x > \zeta(p(w(0)))$.*

Here (i) and (ii) are instances of (iii). For (i), we take $p(x) = Cx^{-\gamma}$ and $w(x) = x^\lambda$ for some $\lambda > 1$, which makes (1.6) equivalent to the statement that $V_1$ does not have a finite exponential moment. For (ii), we take $p(x) = Ce^{-\alpha x}$ and $w(x) = x^{1/\theta}$, which makes (1.6) equivalent to the statement that $E^V((V_1^+)^\theta) = \infty$. In general, (1.6) says that the tails of both $E_1$ and $V_1^+$ are sufficiently fat. The fatter the tails of $E_1$ are, the less fat the tails of $V_1$ need to be.

**2. The deterministic critical point.** Before proving our main results, we need to investigate how a nonrandom potential affects the Markov chain. Many of the basic ideas we need for this already exist, at least when the underlying Markov chain is simple random walk, either rigorously (see, e.g., [13]) or nearly rigorously in the physics literature. But our Markov chain formulation allows an arbitrary excursion-length distribution for the Markov chain, which brings in some complications, so we must go through the details. Let $T_0 = 0$ and for $j \geq 1$, let $T_j$ be the time of the $j$th return to 0 after time 0, if such a return occurs; otherwise, $T_j = \infty$. Let

$$M_E(t) = \langle e^{tE_1} \rangle^X;$$

this is not well defined at $t = 0$ if $P^X(E_1 = \infty) > 0$, so in that case we set $M_E(0) = P^X(E_1 < \infty)$, making $M_E$ left-continuous at 0. Let

$$a_E = a_E(P^X) = \sup\{t \geq 0 : M_E(t) < \infty\}, \qquad a'_E = \lim_{t \nearrow a_E} (\log M_E)'(t).$$

It is easy to see that when $E_1 < \infty$ a.s.,

$$(2.1) \quad a_E = \sup\{s \geq 0 : P^X(E_1 = n) \leq e^{-sn} \text{ for all sufficiently large } n\}.$$

We say that the state 0 is *exponentially recurrent* for $\{X_i\}$ if $a_E > 0$. We write $u_c^d$ for $u_c^d(1)$, which is all we need to consider since $u_c^d(\beta) = \beta^{-1} u_c^d(1)$. We write $m_E$ for the mean of $E_1$, given $E_1 < \infty$. Let

$$J_E(t) = \sup_x (tx - \log M_E(x)).$$

Of course, if $E_1 < \infty$ a.s., then $J_E = I_E$, but they do differ in one case: if $P^X(E_1 = \infty) > 0$ and $P^X(n \leq E_1 < \infty)$ decays exponentially in $n$, then, for all $t > m_E$, we have $-\log P^X(E_1 < \infty) = -\log M_E(0) = J_E(t) < I_E(t)$.

We will show that there exists $C = C(\beta, u)$ such that

$$\lim_n \mu_{[0,n]}^{\beta, u}\left(\frac{L_n}{n} \in (C - \varepsilon, C + \varepsilon)\right) = 1 \qquad \text{for all } \varepsilon > 0;$$

we call $C$ the *contact fraction*. The contact fraction is nondecreasing in $u$, and equal to 0 for $u < u_c^d$. The transition at $u_c^d$ is *first order* if $C(\beta, \cdot)$ is discontinuous at $u_c^d$.

Throughout the paper, $c, c_1, c_2, \ldots$ are unspecified constants; $c$ may take different values at different appearances. Our main result in this section is the following.

THEOREM 2.1. *Let $\{X_i\}$ be an irreducible Markov chain. If 0 is exponentially recurrent for $\{X_i\}$, then*

$$u_c^d = -\log M_E(a_E) \in [-\infty, 0),$$

*the contact fraction exists for all $\beta > 0, u \in \mathbb{R}$ and the transition at $u_c^d$ is not first order. If $0$ is not exponentially recurrent for $\{X_i\}$, then*

$$u_c^d = -\log M_E(a_E) = -\log P^X(E_1 < \infty) \in [0, \infty),$$

*the contact fraction exists for all $\beta > 0$ and $u \neq u_c^d$, and (when $u_c^d$ is finite) the transition is first order if and only if $m_E < \infty$. In both cases, provided $u_c^d$ is finite, there exists $\varphi(\delta) = o(\delta)$ as $\delta \to 0$ such that, for sufficiently small $\delta$, we have*

$$(2.2) \qquad \mu_{[0,n]}^{1,u_c^d}(L_n \leq \delta n) = e^{-o(n)} \qquad \text{as } n \to \infty,$$

*and for large $n$,*

$$(2.3) \qquad \mu_{[0,n]}^{1,u_c^d}(L_n \geq \delta n) \geq e^{-\varphi(\delta)n}.$$

*Finally, the free energy $f^d(\beta, u)$ is given by*

$$(2.4) \qquad \beta f^d(\beta, u) = \lim_n \frac{1}{n} \log Z_{[0,n]}(\beta, u) = \sup_{\delta \in (0,1)} (\beta u \delta - \delta J_E(\delta^{-1})),$$

*which is equal to $-a_E$ for all $u < u_c^d$, and is strictly greater than $-a_E$ for all $u > u_c^d$.*

For the proof we will need the following.

LEMMA 2.2. *Let $\{X_i\}$ be an irreducible Markov chain. For $x \in (0, 1]$, define*

$$g(x) = x J_E(x^{-1}).$$

*For all $0 \leq \delta < \eta \leq 1$, we have*

$$\lim_n \frac{1}{n} \log P^X(\delta n < L_n \leq \eta n) = -\inf_{\delta < x \leq \eta} g(x).$$

The reason that $J_E$ and not $I_E$ appears in the definition of $g$ in this lemma is essentially the following. Suppose $P^X(E_1 = \infty) > 0$, suppose $P^X(n \leq E_1 < \infty)$ decays exponentially in $n$, and suppose we condition on $L_n \geq \delta n$ for some $\delta \in (0, m_E^{-1})$. With high probability, we will not see $\delta n$ excursions of average length near $\delta^{-1}$, which would have a cost per excursion of $I_E(\delta^{-1})$. Instead we will see about $\delta n$ excursions of average length near $m_E^{-1}$ (the last excursion thus ending well before time $n$), followed by an escape to infinity, as this has a lower cost per excursion of $-\log P^X(E_1 < \infty) = J_E(\delta^{-1})$. This is reflected in (2.7) in the proof.

PROOF OF LEMMA 2.2. Let $b = -\log P^X(E_1 < \infty)$. The sup in the expression

$$J_E(t) = \sup_{x < a'_E} (xt - \log M_E(x))$$

occurs for $t < a'_E$ at a unique point $x_0(t) < a_E$ satisfying

$$J'_E(t) = x_0(t), \qquad J_E(t) - tJ'_E(t) = -\log M_E(x_0(t)),$$

$$x_0(t) \to a_E \text{ as } t \to \infty,$$

and occurs uniquely at $x_0(t) = a_E$ for $t \geq a'_E$. In addition, $J_E$ is strictly decreasing and strictly convex on $[0, m_E^{-1}]$, and convex and nondecreasing on $[m_E^{-1}, \infty)$, with $J_E(m_E^{-1}) = b$. Further,

$$(2.5) \qquad \frac{J_E(t)}{t} \to a_E \qquad \text{as } t \to \infty.$$

Define $g(0) = a_E$. We have

$$g'(\delta) = -\log M_E(x_0(\delta^{-1})), \qquad \delta > 0,$$

which is nondecreasing, so $g \geq 0$ is convex, and $g$ is strictly convex and strictly increasing on $[m_E^{-1}, 1]$. Further, $g$ is continuous at 0 by (2.5), and

$$g'(0+) = -\log M_E(a_E).$$

Now let $0 \leq \delta < \eta \leq 1, r = \lfloor \delta n \rfloor + 1, s = \lfloor \eta n \rfloor + 1$, so that

$$(2.6) \qquad P^X(\delta n < L_n \leq \eta n) = P^X(T_r \leq n < T_s).$$

From here we consider the two cases in the theorem statement separately.

*Case* 1. Suppose that 0 is not exponentially recurrent for $\{X_i\}$. We then have

$$a_E = 0, \qquad a'_E = m_E, \qquad \log M_E(a_E) = -b,$$

$$J_E(t) = b \quad \text{for } t \geq m_E, \qquad J_E(t) > b \quad \text{for } t < m_E$$

and

$$g(\delta) = b\delta \qquad \text{for } \delta \in [0, m_E^{-1}].$$

If $m_E = \infty$, then $g'(\delta) > b$ for all $\delta > 0$. For $0 \leq \theta \leq m_E^{-1}$ and $j = \lfloor \theta n \rfloor + 1$, we have

$$(2.7) \quad \lim_n \frac{1}{n} \log P^X(T_j \leq n) = \lim_n \frac{1}{n} \log P^X(E_1 < \infty)^j = -b\theta = -g(\theta),$$

while, for $m_E^{-1} \leq \theta \leq 1$, we have

$$\lim_n \frac{1}{n} \log P^X(T_j \leq n) = \lim_n \frac{1}{n} \log P^X\left(\frac{T_j}{j} \leq \frac{n}{j}\right)$$

$$(2.8) \qquad\qquad = -\theta I_E(\theta^{-1}) = -g(\theta).$$

Hence, if $b > 0$ (so that $g$ is strictly increasing), we have

$$\lim_n \frac{1}{n} \log P^X(\delta n < L_n \leq \eta n) = \lim_n \frac{1}{n} \log P^X\left(\frac{T_r}{r} \leq \frac{n}{r}\right)$$
$$= -g(\delta)$$
$$= -\inf_{\delta < x \leq \eta} g(x).$$

This also holds if $b = 0, \delta \geq m_E^{-1}$, so if $b = 0, \delta \leq m_E^{-1} \leq \eta$, we have

$$
\begin{aligned}
0 &\geq \lim_n \frac{1}{n} \log P^X(\delta n < L_n \leq \eta n) \\
&\geq \lim_n \frac{1}{n} \log P^X(m_E^{-1} n < L_n \leq \eta n) \\
&= -g(m_E^{-1}) \\
&= 0,
\end{aligned}
\tag{2.9}
$$

and then

$$\lim_n \frac{1}{n} \log P^X(\delta n < L_n \leq \eta n) = 0 = -\inf_{\delta < x \leq \eta} g(x). \tag{2.10}$$

Finally, we consider $b = 0, \delta \leq \eta \leq m_E^{-1}$. In that case, since $a_E = 0$, there exists a sequence $q_n \to \infty$ satisfying $q_n = o(n), \log P^X(E_1 = q_n) = o(q_n)$, and we have using (2.10) with $\eta$ replaced by 1 that

$$
\begin{aligned}
0 &\geq \lim_n \frac{1}{n} \log P^X(\delta n < L_n \leq \eta n) \\
&\geq \lim_n \frac{1}{n} \log P^X(T_j \leq n, E_{j+m} = q_n \text{ for all } 1 \leq m \leq n/q_n) \\
&\geq \lim_n \frac{1}{n} \log P^X(T_j \leq n) + \lim_n \frac{1}{q_n} \log P^X(E_1 = q_n) \\
&= \lim_n \frac{1}{n} \log P^X(\delta n < L_n \leq n) \\
&= -\inf_{\delta < x \leq 1} g(x) \\
&= 0,
\end{aligned}
\tag{2.11}
$$

so (2.10) again holds.

*Case* 2. Suppose that 0 is exponentially recurrent for $\{X_i\}$. We then have

$$J_E = I_E, \qquad a_E > 0, \qquad \log M_E(a_E) > 0, \qquad J_E(t) > 0 \qquad \text{for all } t \neq m_E,$$

$$g \text{ strictly decreasing on } [0, m_E] \text{ and } g'(0+) = -\log M_E(a_E).$$

If $\delta < \eta \leq m_E^{-1}$, then, using (2.6) and strict monotonicity of $g$,

$$
\begin{aligned}
(2.12) \qquad & \lim_n \frac{1}{n} P^X(\delta n < L_n \leq \eta n) \\
&= \lim_n \frac{1}{n} \log\left( P^X\left(\frac{T_s}{s} > \frac{n}{s}\right) - P^X\left(\frac{T_r}{r} > \frac{n}{r}\right) \right) \\
&= -g(\eta) \\
&= -\inf_{\delta < x \leq \eta} g(x).
\end{aligned}
$$

Similarly, if $m_E^{-1} \leq \delta < \eta$, then

$$
\begin{aligned}
(2.13) \qquad & \lim_n \frac{1}{n} \log P^X(\delta n < L_n \leq \eta n) \\
&= \lim_n \frac{1}{n} \log\left( P^X\left(\frac{T_r}{r} \leq \frac{n}{r}\right) - P^X\left(\frac{T_s}{s} \leq \frac{n}{s}\right) \right) \\
&= -g(\delta) \\
&= -\inf_{\delta < x \leq \eta} g(x).
\end{aligned}
$$

Finally, if $\delta < m_E^{-1} < \eta$, then $P^X(T_r \leq n < T_s) \to 1$, so, from (2.6),

$$
\lim_n \frac{1}{n} \log P^X(\delta n < L_n \leq \eta n) = 0 = -g(m_E^{-1}) = -\inf_{\delta < x \leq \eta} g(x). \qquad \square
$$

Let $Z_{[0,n]}(\beta, u, \delta-)$ and $Z_{[0,n]}(\beta, u, \delta+)$ denote the contributions to the partition function $Z_{[0,n]}(\beta, u)$ from trajectories with $L_n \leq \delta n$ and $L_n > \delta n$, respectively.

PROOF OF THEOREM 2.1. Fix $u$. It suffices to consider $\beta = 1$, since $Z_{[0,n]}(\beta, u) = Z_{[0,n]}(1, \beta u)$. Let $b$ be as in the proof of Lemma 2.2 and let $f(x) = ux - g(x)$. It follows straightforwardly from Lemma 2.2 that, for $\delta \in (0, 1)$,

$$
\begin{aligned}
(2.14) \qquad & \lim_n \frac{1}{n} \log Z_{[0,n]}(1, u, \delta+) = \sup_{x \in [\delta, 1]} f(x), \\
& \lim_n \frac{1}{n} \log Z_{[0,n]}(1, u, \delta-) = \sup_{x \in [0, \delta]} f(x).
\end{aligned}
$$

If $0$ is not exponentially recurrent, then since $g$ is convex with $g'(0+) = -\log M_E(a_E) = b$ and $g(\delta) = b\delta$ for $0 \leq \delta \leq m_E^{-1}$, we see that $\sup_{x \in [0,1]} f(x)$ occurs uniquely at $x = 0$ if $u < b$ and uniquely at some $x > m_E^{-1}$ if $u > b$, which with (2.14) shows that $u_c^d = b$ and proves (2.2), (2.4), the existence of the contact fraction and the first-order criterion.

If 0 is exponentially recurrent, then since $g$ is strictly convex with $g'(0+) = -\log M_E(a_E)$, we see that $\sup_{x \in [0,1]} f(x)$ occurs uniquely at $x = 0$ if $u < -\log M_E(a_E)$, and uniquely at some $x_u > m_E^{-1}$ if $u > -\log M_E(a_E)$, with $x_u \to 0$ as $u \searrow -\log M_E(a_E)$. With (2.14), this shows that $u_c^d = -\log M_E(a_E)$, proves (2.2), (2.4) and the existence of the contact fraction, and (when $u_c^d$ is finite) shows that the transition is not first order.

In both cases, at $u = u_c^d$, $\sup_{x \in [0,1]} f(x)$ occurs, not necessarily uniquely, at $x = 0$. We have $f'(0+) = 0$, so that $f(0) - \sup_{x \in [\delta, 1]} f(x) = f(0) - f(\delta) = o(\delta)$ as $\delta \to 0$. With (2.14), this proves (2.3).  $\square$

For $u < 0$, one can interpret $1 - e^u$ as a probability for the Markov chain to be killed, each time it visits state 0; $\mu_{[0,n]}^{1,u}$ then becomes the distribution of the chain conditioned on its still being alive at time $n$. For $u > 0$, the killing instead occurs at all states other than 0, with probability $1 - e^{-u}$. This turns the question of pinning into a question about quasistationary distributions, which makes Theorem 2.1 closely related to results in [21].

Let $r_1 < r_2$ be the two smallest values in the set $\{k \geq 1 : P^X(E_1 = k) > 0\}$ of possible first-return times. We refer to $(T_{2j-2}, T_{2j}]$ as the $j$th *block*. We say the $j$th block is *good* if $T_{2j} - T_{2j-2} = r_1 + r_2$, and *bad* otherwise, and define

$$p_g = P^X(\text{block 1 is good}) = 2P^X(E_1 = r_1)P^X(E_1 = r_2),$$

$$G_n = |\{j \geq 1 : T_{2j} \leq n, \text{ block } j \text{ is good}\}|,$$

$$\delta_E = \frac{1}{\langle E_1 \rangle^X}.$$

$\delta_E$ represents the "natural frequency" for returns to 0.

In proving our results for random potentials, we will be interested in probabilities

$$\mu_{[0,n]}^{\beta,u}(G_n \geq c_1 p_g \delta n | L_n = k_n).$$

We observe that this probability does not depend on $\beta$ or $u$, so we need only consider $u = 0$, which makes $\mu_{[0,n]}^{\beta,u} = P_{[0,n]}^X$. In the exponentially recurrent case, if $\delta > 0$, and $k_n \geq \delta n$ is much smaller than $\delta_E n$, the chain conditioned on $L_n = k_n$ may make a large number of unusually long excursions. The question is, under such conditioning, could the chain also then typically have an unusually small proportion of short excursions, and, consequently, $G_n$ be typically much smaller than $p_g \delta n$? The next lemma shows that the answer is no, when $a_E < \infty$ and $M_E(a_E) < \infty$.

LEMMA 2.3.    *Let $\{X_i\}$ be an irreducible Markov chain with $a_E < \infty$ and $M_E(a_E) < \infty$. There exists $c_1 > 0$ such that, for all $\delta > 0$ and all sequences*

$\{k_n\}$ with $\delta n \le k_n \le (1-\delta) r_1^{-1} n$, we have

$$P^X(G_n \ge c_1 p_g \delta n | L_n = k_n) \to 1 \qquad \text{as } n \to \infty.$$

PROOF. Fix $\delta$ and $\{k_n\}$ as in the lemma statement. The returns to 0 form the arrivals of a renewal process with interarrival times $E_i$, so we may reformulate the lemma as a statement about such a process. Thus, $T_k$ denotes the time of the $k$th renewal, and we let $\tilde{G}_n$ denote the number of good blocks among the first $\lfloor \delta n/2 \rfloor$ blocks, so that, given $L_n = k_n$, we have $G_n \ge \tilde{G}_n$. It follows easily from basic large deviation theory that, for $r_1 < b < \min(r_1/(1-\delta), m_E)$, we have

$$P^X(T_{k_n} \le b k_n | L_n = k_n) \to 0 \qquad \text{as } n \to \infty.$$

Fix such a $b$ and let $y_n \to \infty$ with $b k_n < y_n \le n$; it is thus sufficient to show that, for all such $\{y_n\}$, we have

$$(2.15) \qquad P^X(\tilde{G}_n \ge c_1 p_g \delta n | T_{k_n} = y_n) \to 1 \qquad \text{as } n \to \infty.$$

We may assume that $E_1 < \infty$ a.s. We tilt the distribution of $E_1$, defining the measures

$$Q^t(\cdot) = \frac{\langle e^{tE_1} \delta_{\{E_1 \in \cdot\}} \rangle^X}{\langle e^{tE_1} \rangle^X}$$

whenever $M_E(t) < \infty$, and let $\nu^t$ denote the distribution of the renewal process with interarrival distribution $Q^t$. We observe that the probability (2.15) is unchanged if we replace $P^X$ (or, equivalently, $\nu^0$) with $\nu^t$, for arbitrary $t$ satisfying $M_E(t) < \infty$. By considering subsequences, we may assume that $y_n/k_n \to \kappa$ for some $\kappa \in [b, \delta^{-1}]$.

Since $(\log M_E)'(0-) = m_E > b > r_1 = \lim_{t \to -\infty} (\log m_E)'(t)$, we can define $b_E$ by

$$(\log M_E)'(b_E) = b,$$

and let $t \in (b_E, a_E)$. We have

$$(2.16) \qquad \nu^t(\tilde{G}_n \le c_1 p_g k_n | T_{k_n} = y_n) \le \frac{\nu^t(\tilde{G}_n \le c_1 p_g k_n)}{\nu^t(T_{k_n} = y_n)}$$

and for $i = 1, 2$,

$$Q^t(r_i) = \frac{e^{tr_i} P^X(E_1 = r_i)}{M_E(t)} \ge \frac{e^{b_E r_i} P^X(E_1 = r_i)}{M_E(a_E)},$$

so

$$\nu^t(\text{block } j \text{ is good}) = 2Q^t(r_1) Q^t(r_2)$$
$$\ge 2\left(\frac{e^{b_E r_1}}{M_E(a_E)}\right)^2 P^X(E_1 = r_1) P^X(E_1 = r_2)$$
$$= \left(\frac{e^{b_E r_1}}{M_E(a_E)}\right)^2 p_g.$$

Thus, taking $0 < c_1 < \frac{1}{2}(e^{b_E r_1}/M_E(a_E))^2$, we get that, for some $\gamma > 0$ depending only on $P^X$ and $c$,

$$(2.17) \qquad \nu^t(\tilde{G}_n \leq c_1 p_g \delta n) \leq e^{-\gamma \delta n} \qquad \text{for all } t \in (b_E, a_E).$$

We now need to choose a $t = t_n$ that gives a good lower bound on the denominator in (2.16). Consider first the case $\kappa < a'_E$. For large $n$, there exists $t_n \in (b_E, a_E)$ with $(\log M_E)'(t_n) = y_n/k_n$, and, hence, under $\nu^{t_n}$, the $E_i$ have mean $y_n/k_n$. There also exists $t_\infty \in [b_E, a_E)$ with $(\log M_E)'(t_\infty) = \kappa$, and all moments (including exponential) of $Q^{t_n}$ converge to those of $Q^{t_\infty}$. The standard proof of the local central limit theorem (see, e.g., [7]) carries over to this situation and shows that

$$(2.18) \qquad n^{1/2} \nu^{t_n}(T_{k_n} = y_n) \to \frac{1}{\sqrt{2\pi} \sigma_\infty} \qquad \text{as } n \to \infty,$$

where $\sigma_\infty$ is the standard deviation of the measure $Q^{t_\infty}$. With (2.16) and (2.17), this proves (2.15).

Next we consider $\kappa \geq a'_E$. Here we cannot necessarily tilt the distribution of $E_1$ to change the mean to $y_n/k_n$. Instead we tilt to obtain a lower mean, and force the average excursion length up to $y_n/k_n$ using a small number of much longer excursions. Specifically, let $0 < \varepsilon < \delta \gamma / 8$, let $m_t$ be the mean of $Q^t$, and take $t = a_E - \varepsilon$. Then $m_t < a'_E \leq \kappa$ so $m_t < y_n/k_n$ for all sufficiently large $n$. By (2.1), we have

$$0 = \sup\{s \geq 0 \colon Q^{a_E}(E_1 = n) \leq e^{-sn} \text{ for all sufficiently large } n\},$$

so we can choose $q > \kappa$ satisfying

$$Q^{a_E}(E_1 = q) \geq e^{-\varepsilon q},$$

and we then have

$$(2.19) \qquad Q^t(E_1 = q) \geq e^{-2\varepsilon q}.$$

Now let

$$j_n = \min\{j \geq 0 \colon j m_t + (k_n - j)q \leq y_n\},$$

so that, for large $n$, for some $0 \leq l_n < q - m_t$,

$$0 < j_n < k_n, \qquad k_n - j_n \leq \frac{n}{q} \quad \text{and} \quad y_n = j_n m_t + (k_n - j_n)q + l_n,$$

with $j_n \to \infty$. Then using the local CLT again along with (2.19), for large $n$,

$$
\begin{aligned}
(2.20) \qquad \nu^t(T_{k_n} = y_n) &\geq \nu^t(T_{j_n} = j_n \tilde{m}_t + l_n)\nu^t(T_{k_n - j_n} = q(k_n - j_n)) \\
&\geq \nu^t(T_{j_n} = j_n \tilde{m}_t + l_n)Q^t(E_1 = q)^{k_n - j_n} \\
&\geq \frac{e^{-(\gamma \delta j_n/4 + 2\varepsilon q(k_n - j_n))}}{2\sqrt{2\pi} \sigma_t} \\
&\geq e^{-3\gamma \delta n/4},
\end{aligned}
$$

where $\sigma_t$ denotes the standard deviation of $Q^t$. With (2.17) and (2.16), this proves (2.15). □

**3. Self-averaging of the free energy and critical point.** In this section we establish the existence of a well-defined nonrandom quenched critical point. Self-averaging of the free energy is established for other polymer models in [18] and [23]. Let $Z_{[0,n]}^{\{V_i\}}(\beta, u, \delta-)$ and $Z_{[0,n]}^{\{V_i\}}(\beta, u, \delta+)$ denote the contributions to the partition function $Z_{[0,n]}^{\{V_i\}}(\beta, u)$ from trajectories with $L_n \leq \delta n$ and $L_n > \delta n$, respectively. As before, we omit the $\{V_i\}$ when $V_i \equiv 0$.

THEOREM 3.1. *Let $\{X_i\}$ be an irreducible Markov chain with state space containing a state 0, and let the potential $\{V_i\}$ at 0 be i.i.d. with mean 0. Then for $\beta > 0, u \in \mathbb{R}$, there exists a nonrandom constant $f^q(\beta, u)$ such that*

$$\lim_n \frac{1}{n} \log Z_{[0,n]}^{\{V_i\}}(\beta, u) = \beta f^q(\beta, u), \qquad P^V\text{-a.s.}$$

*and for $\beta > 0$, there exists a constant $u_c^q(\beta)$ such that $u_c(\beta, \{V_i\}) = u_c^q(\beta)$, $P^V$-a.s. and*

$$\beta f^q(\beta, u) = -a_E \qquad \text{for all } u < u_c^q(\beta),$$
$$\beta f^q(\beta, u) > -a_E \qquad \text{for all } u > u_c^q(\beta).$$

PROOF. Fix $\beta > 0, u \in \mathbb{R}$ and define the random variables $F_\pm(\beta, u)$ by

$$\beta F_+(\beta, u) = \limsup_n \frac{1}{n} \log Z_{[0,n]}^{\{V_i\}}(\beta, u),$$

$$\beta F_-(\beta, u) = \liminf_n \frac{1}{n} \log Z_{[0,n]}^{\{V_i\}}(\beta, u),$$

$$U_0 = \inf\left\{u \in \mathbb{R} : \limsup_n \mu_{[0,n]}^{\beta, u, \{V_i\}}(L_n \geq \delta n) > 0 \text{ for some } \delta > 0\right\}.$$

It is easy to see that $F_\pm(\beta, u)$ and $U_0$ are tail random variables of the sequence $\{V_i\}$, so there exist $f_\pm^q(\beta, u)$ and $u_0$ such that $F_\pm(\beta, u) = f_\pm^q(\beta, u)$ a.s. and $U_0 = u_0$ a.s. Fix $M > 0$ to be specified and consider the truncated potential $\tilde{V}_i = (V_i \wedge M) \vee (-M)$. From (2.14) in the proof of Theorem 2.1, we have

$$-a_E = -g'(0) = \lim_{\delta \to 0} \lim_n \frac{1}{n} \log Z_{[0,n]}(\beta, u, \delta-).$$

Observe that, for all $u \in \mathbb{R}, \delta > 0$,

$$|\log Z_{[0,n]}^{\{\tilde{V}_i\}}(\beta, u, \delta-) - \log Z_{[0,n]}^{\{V_i\}}(\beta, u, \delta-)| \leq \beta \sum_{i=1}^n |V_i| \delta_{\{|V_i| > M\}}$$

and hence

$$
(3.1) \quad \limsup_n \left| \frac{1}{n} \log Z_{[0,n]}^{\{\tilde{V}_i\}}(\beta, u, \delta-) - \frac{1}{n} \log Z_{[0,n]}^{\{V_i\}}(\beta, u, \delta-) \right|
$$
$$
\leq \beta \langle |V_1| \delta_{\{|V_1|>M\}} \rangle^V, \qquad P^V\text{-a.s.}
$$

Therefore,

$$
\begin{aligned}
(3.2) \quad \beta f_-^q(\beta, u) &\geq \liminf_n \frac{1}{n} \log Z_{[0,n]}^{\{V_i\}}(\beta, u, \delta-) \\
&\geq \liminf_n \frac{1}{n} \log Z_{[0,n]}^{\{\tilde{V}_i\}}(\beta, u, \delta-) - \beta \langle |V_1| \delta_{\{|V_1|>M\}} \rangle^V \\
&\geq \lim_n \frac{1}{n} \log Z_{[0,n]}(\beta, u, \delta-) - \delta(|u| + M) - \beta \langle |V_1| \delta_{\{|V_1|>M\}} \rangle^V \\
&\geq -a_E - \delta(|u| + M) - \beta \langle |V_1| \delta_{\{|V_1|>M\}} \rangle^V.
\end{aligned}
$$

We can take $M$ large and then $\delta$ small, so we get $\beta f_-^q(\beta, u) \geq -a_E$. In the other direction, for $u < u_0$ and $\delta > 0$, we have

$$
(3.3) \quad \beta f_+^q(\beta, u) = \limsup_n \frac{1}{n} \log Z_{[0,n]}^{\{V_i\}}(\beta, u, \delta-),
$$

so we obtain similarly to (3.2) that $\beta f_+^q(\beta, u) \leq -a_E$, and therefore,

$$
(3.4) \quad \lim_n \frac{1}{n} \log Z_{[0,n]}^{\{V_i\}}(\beta, u) = -a_E \qquad \text{a.s. for all } u < u_0.
$$

It remains to consider $u > u_0$. Defining

$$
\Delta_0(u) = \sup \left\{ \delta > 0 : \limsup_n \mu_{[0,n]}^{\beta, u, \{V_i\}}(L_n \geq \delta n) > 0 \right\},
$$

we see as with $U_0$ that there exists $\delta_0(u) > 0$ such that $\Delta_0(u) = \delta_0(u)$ a.s., and $\delta_0(u)$ is an increasing function of $u$. Fix $u_0 < v < u$ and $0 < \eta < \delta < \delta_0(v)$. Then

$$
\limsup_n \left( \frac{1}{n} \log Z_{[0,n]}^{\{V_i\}}(\beta, v, \delta+) - \frac{1}{n} \log Z_{[0,n]}^{\{V_i\}}(\beta, v, \eta-) \right) \geq 0,
$$

so, provided $\eta$ and $v - u_0$ are small enough,

$$
\begin{aligned}
(3.5) \quad &\limsup_n \left( \frac{1}{n} \log Z_{[0,n]}^{\{V_i\}}(\beta, u, \delta+) - \frac{1}{n} \log Z_{[0,n]}^{\{V_i\}}(\beta, u, \eta-) \right) \\
&\qquad \geq (u - v)(\delta - \eta) \\
&\qquad > \frac{3}{4}(u - u_0)\delta.
\end{aligned}
$$

Arguing as in (3.2), we see that, for $M$ sufficiently large and $\eta$ sufficiently small, it follows from (3.5) that

$$(3.6) \qquad \limsup_n \left( \frac{1}{n} \log Z_{[0,n]}^{\{\tilde{V}_i\}}(\beta, u, \delta+) + a_E \right) > \frac{1}{2}(u-u_0)\delta.$$

We would like to use superadditivity of the mean $\langle \log Z_{[0,n]}^{\{\tilde{V}_i\}}(\beta, u, \delta+) \rangle^V$ to help us conclude that the lim sup in (3.6) is actually a limit, but unfortunately this sequence is not obviously superadditive unless we restrict to paths which end at state 0 at time $n$. To circumvent this difficulty, we proceed as follows. Let $Z_{[0,n]}^{\{\tilde{V}_i\}}(\beta, u, \delta+, k)$ denote the contribution to $Z_{[0,n]}^{\{\tilde{V}_i\}}(\beta, u, \delta+)$ from trajectories with the last visit to 0 in $[0,n]$ at time $k$. There exists $k_n \geq \delta n$ (depending on $\{V_i\}$) such that, letting $G_E$ denote the distribution function of $E_1$ under $P^X$,

$$(3.7) \qquad \begin{aligned} Z_{[0,k_n]}^{\{\tilde{V}_i\}}(\beta, u, \delta+, k_n)(1 - G_E(n-k_n)) &= Z_{[0,n]}^{\{\tilde{V}_i\}}(\beta, u, \delta+, k_n) \\ &\geq \frac{1}{n} Z_{[0,n]}^{\{\tilde{V}_i\}}(\beta, u, \delta+). \end{aligned}$$

From (2.1) and (3.7), we obtain

$$(3.8) \qquad \begin{aligned} &\frac{1}{n} \log Z_{[0,n]}^{\{\tilde{V}_i\}}(\beta, u, \delta+) \\ &\leq \frac{1}{n} \log Z_{[0,k_n]}^{\{\tilde{V}_i\}}(\beta, u, \delta+, k_n) + \frac{1}{n} \log(1 - G_E(n-k_n)) + \frac{\log n}{n} \\ &\leq \frac{1}{n} \log Z_{[0,k_n]}^{\{\tilde{V}_i\}}(\beta, u, \delta+, k_n) - \left(1 - \frac{k_n}{n}\right) a_E + \frac{1}{8}(u-u_0)\delta^2. \end{aligned}$$

Combining (3.6), (3.7) and (3.8), we obtain

$$(3.9) \qquad \limsup_n \frac{1}{k_n} \log Z_{[0,k_n]}^{\{\tilde{V}_i\}}(\beta, u, \delta+, k_n) > -a_E + \frac{1}{4}(u-u_0)\delta.$$

Analogously to [14], from Azuma's inequality [2], we get that, for some $K > 0$ depending on $M$, for all $\delta > 0, k \geq 1$,

$$\begin{aligned} P^V(|\log &Z_{[0,k]}^{\{\tilde{V}_i\}}(\beta, u, \delta+, k) \\ &- \langle \log Z_{[0,k]}^{\{\tilde{V}_i\}}(\beta, u, \delta+, k) \rangle^V| > K k^{1/2} \log k) \leq k^{-2}. \end{aligned}$$

With the Borel–Cantelli lemma and (3.9), this shows that, for some deterministic $\{k_n\}$,

$$(3.10) \qquad \limsup_n \frac{1}{k_n} \langle \log Z_{[0,k_n]}^{\{\tilde{V}_i\}}(\beta, u, \delta+, k_n) \rangle^V > -a_E + \frac{1}{4}(u-u_0)\delta.$$

Therefore, we can choose a fixed $m$ satisfying

$$(3.11) \qquad \frac{1}{m} \langle \log Z_{[0,m]}^{\{\tilde{V}_i\}}(\beta, u, \delta+, m) \rangle^V > -a_E + \frac{1}{4}(u - u_0)\delta.$$

Observe that the sequence

$$b_j = \langle \log Z_{[0,jm]}^{\{\tilde{V}_i\}}(\beta, u, \delta+, jm) \rangle^V, \qquad j \geq 1,$$

is superadditive, because, for $j, k \geq 1$,

$$Z_{[0,(j+k)m]}^{\{\tilde{V}_i\}}(\beta, u, \delta+, (j+k)m)$$

$$\geq Z_{[0,jm]}^{\{\tilde{V}_i\}}(\beta, u, \delta+, jm) Z_{[jm,(j+k)m]}^{\{\tilde{V}_i\}}(\beta, u, \delta+, (j+k)m),$$

where the last partition function is for trajectories $\{x_i, jm \leq i \leq (j+k)m\}$ with $x_{jm} = x_{(j+k)m} = 0$ and at least $\delta n$ returns to 0. Therefore, the limit

$$f^q(\beta, u, \delta, M) = \frac{1}{\beta} \lim_{j \to \infty} \frac{b_j}{jm}$$

exists, and $b_j/jm \leq \beta f^q(\beta, u, \delta, M)$ for all $j$, so

$$(3.12) \qquad \beta f^q(\beta, u, \delta, M) > -a_E + \frac{1}{4}(u - u_0)\delta.$$

It follows easily from boundedness of $\tilde{V}_1$ that in fact the convergence occurs for the full sequence:

$$(3.13) \qquad \lim_n \frac{1}{n} \langle \log Z_{[0,n]}^{\{\tilde{V}_i\}}(\beta, u, \delta+, n) \rangle^V = \beta f^q(\beta, u, \delta, M).$$

For every choice of $\delta n \leq k_n \leq n$, we have by (2.1), the equality in (3.7) and (3.12) that

$$(3.14) \quad \begin{aligned} \limsup_n &\frac{1}{n} \langle \log Z_{[0,n]}^{\{\tilde{V}_i\}}(\beta, u, \delta+, k_n) \rangle^V \\ &\leq \limsup_n \Big( \frac{k_n}{n} \frac{1}{k_n} \langle \log Z_{[0,k_n]}^{\{\tilde{V}_i\}}(\beta, u, \delta+, k_n) \rangle^V \\ &\qquad\qquad + \Big(1 - \frac{k_n}{n}\Big) \frac{1}{n - k_n} \log P^X(E_1 > n - k_n) \Big) \\ &\leq \beta f^q(\beta, u, \delta, M), \end{aligned}$$

which with (3.13) and the inequality in (3.7) shows that

$$(3.15) \qquad \lim_n \frac{1}{n} \langle \log Z_{[0,n]}^{\{\tilde{V}_i\}}(\beta, u, \delta+) \rangle^V = \beta f^q(\beta, u, \delta, M),$$

for all $\delta < \delta_0(u)$. With $\delta$ fixed, taking $M$ sufficiently large and then $\eta$ sufficiently small, we obtain as in (3.2), using also (3.12), that

$$
\begin{aligned}
(3.16) \quad & \limsup_n \frac{1}{n} \langle \log Z_{[0,n]}^{\{\tilde V_i\}}(\beta, u, \eta-) \rangle^V \\
& \leq -a_E + \eta(|u| + M) + \beta \langle |V_1| \delta_{\{|V_1| > M\}} \rangle^V \\
& < \beta f^q(\beta, u, \delta, M) \leq \beta f^q(\beta, u, \eta, M),
\end{aligned}
$$

so that, using (3.16), (3.15) and (3.12),

$$
\begin{aligned}
(3.17) \quad & \lim_n \frac{1}{n} \langle \log Z_{[0,n]}^{\{\tilde V_i\}}(\beta, u) \rangle^V = \lim_n \frac{1}{n} \langle \log Z_{[0,n]}^{\{\tilde V_i\}}(\beta, u, \eta+) \rangle^V \\
& = \beta f^q(\beta, u, \eta, M) \\
& > -a_E + \frac{1}{4}(u - u_0)\delta,
\end{aligned}
$$

meaning that $f^q(\beta, u, \eta, M)$ does not depend on (small) $\eta$, so we denote it $f^q(\beta, u, M)$. Using Azuma's inequality again and (3.17), we get

$$
\lim_n \frac{1}{n} \log Z_{[0,n]}^{\{\tilde V_i\}}(\beta, u) = \beta f^q(\beta, u, M) > -a_E + \frac{1}{4}(u - u_0)\delta \qquad \text{a.s.}
$$

for all $M$, and then from (3.1) with $\delta = 1$, there exists

$$
f^q(\beta, u) = \lim_{M \to \infty} f^q(\beta, u, M)
$$

satisfying

$$
(3.18) \quad \lim_n \frac{1}{n} \log Z_{[0,n]}^{\{V_i\}}(\beta, u) = \beta f^q(\beta, u) \geq -a_E + \frac{1}{4}(u - u_0)\delta \qquad \text{a.s.,}
$$

which establishes the existence of the quenched free energy.

From the definition of $u_0$, there is no pinning at $(\beta, u)$ for $u < u_0$. For $u > u_0$, similarly to (3.16), we have for sufficiently small $\eta$ that

$$
(3.19) \quad \limsup_n \frac{1}{n} \langle \log Z_{[0,n]}^{\{V_i\}}(\beta, u, \eta-) \rangle^V \leq -a_E + \frac{1}{8}(u - u_0)\delta,
$$

which with another application of Azuma's inequality yields

$$
\limsup_n \frac{1}{n} \log Z_{[0,n]}^{\{V_i\}}(\beta, u, \eta-) \leq -a_E + \frac{1}{8}(u - u_0)\delta \qquad \text{a.s.}
$$

This and (3.18) show that for $u > u_0$ there is pinning at $(\beta, u)$. It follows that $u_c(\beta, \{V_i\}) = u_0$ a.s. $\quad \square$

**4. Proof of Theorem 1.5.**   We may assume the Markov chain is aperiodic. We establish pinning by finding $0 < \theta < \delta$ and a set $\Xi$ of pinned trajectories (more precisely, $\Xi \subset \{L_n \geq \delta n\}$) with Gibbs weight exponentially greater than the set $\{L_n \leq \theta n\}$, which, roughly speaking, includes all unpinned trajectories.

As a shorthand we refer to the potentials $V_i$ as *rewards* and say that a reward $V_j$ is *received* by a trajectory $\{x_i\}$ if $x_j = 0$.

We begin by introducing some independence into the sampling done by the Markov chain from the set $\{V_i\}$ by virtue of the times of its returns to 0. Recall that $r_1 < r_2$ are the two smallest values in the set $\{k \geq 1 : P^X(E_1 = k) > 0\}$ of possible first-return times, and $T_j$ is the time of the $j$th return to 0 after time 0, if such a return occurs. Let $W_j^* = T_{2j-2} + r_1$. We refer to $W_j^*$ as the $j$th *target*. When the $j$th block is good we say the $j$th target is *hit* if $T_{2j-1} = W_j^*$, and *missed* otherwise; a missed target means that $T_{2j-1} = T_{2j-2} + r_2$. Given that a block is good, the target is equally likely to be hit or missed, so

$$p_h = P^X(j\text{th target is hit}|\text{block } j \text{ is good}) = 1/2.$$

We use the notation $p_h$ so that the reader may distinguish this probability from other numerical factors that appear in the proof. Conditionally on $\{T_{2j}, j \geq 1\}$, the targets are hit or missed independently.

Fix $\delta > 0$ to be specified, fix $n$, let $c_1$ be as in Lemma 2.3, let

$$J^* = \lfloor c_1 p_g \delta n \rfloor$$

and define the event

$$\Xi = \{x_{[0,n]} \in \Sigma^{[0,n]} : L_n \geq \delta n, G_n \geq J^*\}.$$

We may assume $c_1 \leq 1/4$. For trajectories $\{x_i\} \in \Xi$, we can define $N = N(\{x_i\})$ by stating that the $N$th block is the $J^*$th good block. We then define

$$\mathfrak{G}^* = \mathfrak{G}^*(\{x_i\}) = \{j \leq N : \text{block } j \text{ is good}\},$$

$$\mathfrak{G} = \mathfrak{G}(\{x_i\}) = \{j \leq N : \text{block } j \text{ is good and target } W_j^* \text{ is hit}\},$$

$$\mathcal{R} = \mathcal{R}(\{x_i\}) = \{i \in [0,n] : x_i = 0\},$$

so that $|\mathfrak{G}^*| = J^*$, and define the random sequences

$$\mathcal{W}^* = (W_j^* : j \in \mathfrak{G}^*), \qquad \mathcal{W} = (T_{2j-1} : j \in \mathfrak{G}), \qquad \mathcal{U} = \mathcal{R} \backslash \mathcal{W} = \mathcal{R} \backslash \mathcal{W}^*.$$

For $R \subset [0,n]$, we set

$$S_R = S_R(\{V_i\}) = \sum_{i \in R} V_i.$$

Let $\varepsilon, \alpha > 0$ to be specified. The idea is to condition on the event $\Xi$, on $T_2, T_4, \ldots, T_{2J^*}$ and on the disorder $\{V_i\}$; this makes

$$S_{\mathcal{W}} = \sum_{j \in \mathfrak{S}^*} V_{W_j^*} \delta_{\{W_j^* \text{ is hit}\}}$$

into an i.i.d. sum, and we consider large deviations for this sum in which the average value $S_{\mathcal{W}}/|\mathcal{W}|$ is of order $\varepsilon$. We will need to ensure that, for typical disorders, the remaining rewards $S_{\mathcal{R}} - S_{\mathcal{W}}$ received by the chain at "nontarget" returns to 0 are unlikely to cancel out a large-deviation value of $S_{\mathcal{W}}$. More precisely, a large deviation for $S_{\mathcal{W}}$ of size $\varepsilon p_g p_h \delta n / 4$ needs to imply with high probability a (slightly smaller) large deviation for $S_{\mathcal{R}}$. We also need to ensure that the large-deviation rate function for $S_{\mathcal{W}}$ under the above conditioning is (for typical disorders) not too different from $I_V$. Under that conditioning, the log moment generating function of $V_{W_j^*} \delta_{\{W_j^* \text{ is hit}\}}$ is $\ell(t) = \log(1 + p_h(e^{tV_{W_j^*}} - 1))$, and we define the corresponding mean

$$\psi(t) = \langle \log(1 + p_h(e^{tV_1} - 1)) \rangle^V$$

and an analog of the rate function:

$$\tilde{I}(t) = -\inf_x (\psi(x) - xt).$$

Let $\eta > 0$ to be specified, let $u = u_c^d - \eta$ and define the product measure

$$\mathbb{P} = P_{[0,n]}^V \times \mu_{[0,n]}^{\beta, u} \qquad \text{on } \mathbb{R}^{[0,n]} \times \Sigma^{[0,n]}.$$

Let

$$\tilde{A} = \left\{ (W, v_W) : W \subset [0,n], v_W \in \mathbb{R}^W, \left| \sum_{i \in W} v_i - \varepsilon c_1 p_g p_h \delta n \right| < \alpha n \right\},$$

$$A = \{ (v_{[0,n]}, x_{[0,n]}) \in \mathbb{R}^{[0,n]} \times \Xi : (\mathcal{W}, V_{\mathcal{W}}) \in \tilde{A} \},$$

$$\tilde{B} = \left\{ (W^*, v_{W^*}) : \frac{1}{J^*} \log \mathbb{P}(A | \mathcal{W}^* = W^*, V_{\mathcal{W}^*} = v_{W^*}) \right.$$

$$\left. \in [-\tilde{I}(p_h \varepsilon) - \alpha, -\tilde{I}(p_h \varepsilon) + \alpha], S_{\mathcal{W}^*} \in [-\alpha J^*, \alpha J^*] \right\},$$

$$B = \{ (v_{[0,n]}, x_{[0,n]}) \in \mathbb{R}^{[0,n]} \times \Xi : (\mathcal{W}^*, V_{\mathcal{W}^*}) \in \tilde{B} \},$$

$$\lambda = \frac{\varepsilon c_1 p_g p_h \delta}{4},$$

$$G = \left\{ (v_{[0,n]}, x_{[0,n]}) \in \mathbb{R}^{[0,n]} \times \Xi : \sum_{i \in \mathcal{U}} v_i \geq -\lambda n \right\}.$$

Loosely speaking, $B$ is the event that the rewards at the targets have a typical degree of conduciveness to a large deviation of order $\varepsilon$, $A$ is the event that such a large deviation actually occurs, and $G$ is the event that this large deviation is not canceled out by the rewards received at nontarget locations. We claim that there exists $\nu > 0$ such that

$$(4.1) \qquad \mathbb{P}(G^c | A \cap B) \le e^{-2\nu n}.$$

It suffices to show that, for every $(W^*, v_{W^*}) \in \tilde{A}$, for every $U \subset [0, n]$ such that $U \cap W^* = \phi$, for every $W \subset W^*$ such that $\mathcal{R} = W \cup U$ implies $X_{[0,n]} \in \Xi$, we have

$$(4.2) \qquad \mathbb{P}(G^c | B, \mathcal{W}^* = W^*, \mathcal{W} = W, V_W = v_W, \mathcal{U} = U) \le e^{-2\nu n}.$$

The only part of the conditioning in (4.2) that is relevant to $G^c$ is $\mathcal{U} = U, \mathcal{W} = W$ with $U, W$ disjoint, which ensures that, conditionally, $S_U$ is just an i.i.d. sum of $|U|$ unconditioned variables $V_i$. More precisely, for $U$ as above, we have $\frac{1}{2}\delta n \le |U| \le n$, and the probability in (4.2) is

$$P_{[0,n]}^V(S_U < -\lambda n) \le e^{-2\nu n},$$

which proves (4.1).

Let $Y_n(\theta)$ denote the sum of the $\lfloor \theta n \rfloor$ largest values among $|V_1|, \ldots, |V_n|$. Since $V_1$ has a finite exponential moment, there exist $a(\theta) \searrow 0$ as $\theta \searrow 0$ and $q(\theta) > 0$ for all $\theta > 0$ such that

$$(4.3) \qquad P^V(Y_n(\theta) > a(\theta)n) \le e^{-q(\theta)n} \qquad \text{for all sufficiently large } n.$$

Define

$$g(v_{[0,n]}) = \mathbb{P}(G^c \cap A \cap B | X_{[0,n]} \in \Xi, V_{[0,n]} = v_{[0,n]}),$$

$$h_\theta(v_{[0,n]}) = \langle e^{\beta S_{\mathcal{R}}(\{v_i\})} \delta_{\{L_n \le \theta n\}} \rangle_{[0,n]}^{\beta, u},$$

$$Q_1 = \{v_{[0,n]} \in \mathbb{R}^{[0,n]} : g(v_{[0,n]}) \le e^{-\nu n - (\bar{I}(p_h \varepsilon) - \alpha)J^*}\},$$

$$Q_2 = \{v_{[0,n]} \in \mathbb{R}^{[0,n]} : \mathbb{P}(B | X_{[0,n]} \in \Xi, V_{[0,n]} = v_{[0,n]}) \ge \tfrac{1}{2}\},$$

$$Q_3 = Q_3(\theta) = \{v_{[0,n]} \in \mathbb{R}^{[0,n]} : h_\theta(v_{[0,n]}) \le e^{\beta a(\theta)n}\},$$

where $\theta$ is to be specified. We need to show that $P_{[0,n]}^X(Q_1 \cap Q_2 \cap Q_3)$ is close to 1, as disorder realizations $\{v_i\} \in Q_1 \cap Q_2 \cap Q_3$ are to be considered "good." We will consider $Q_2$ later. For $Q_1$, from the definition of $B$, we have

$$(4.4) \qquad e^{-(\bar{I}(p_h \varepsilon) + \alpha)J^*} \le \mathbb{P}(A | B, V_{[0,n]} = v_{[0,n]}) \qquad \text{for all } v_{[0,n]},$$

and

$$(4.5) \qquad \mathbb{P}(A | B) \le e^{-(\bar{I}(p_h \varepsilon) - \alpha)J^*},$$

so, using (4.1),

$$(4.6) \qquad \langle g(V_{[0,n]}) \rangle_{[0,n]}^V = \mathbb{P}(G^c \cap A \cap B | X_{[0,n]} \in \Xi) \leq e^{-2\nu n - (\tilde{I}(p_h \varepsilon) - \alpha) J^*}.$$

Therefore,

$$(4.7) \qquad P_{[0,n]}^V(Q_1^c) \leq e^{\nu n + (\tilde{I}(p_h \varepsilon) - \alpha) J^*} \langle g(V_{[0,n]}) \rangle_{[0,n]}^V \leq e^{-\nu n}.$$

For $Q_3$, we have $h_\theta(V_{[0,n]}) \leq e^{\beta Y_n(\theta)}$ so, by (4.3),

$$(4.8) \qquad P_{[0,n]}^V(Q_3^c) \leq e^{-q(\theta)n}.$$

For $v_{[0,n]} \in Q_2$, we have by (4.4) that

$$(4.9) \qquad \mathbb{P}(A \cap B | X_{[0,n]} \in \Xi, V_{[0,n]} = v_{[0,n]}) \geq \tfrac{1}{2} e^{-(\tilde{I}(p_h \varepsilon) + \alpha) J^*}.$$

Provided $\alpha < \nu/2$, for $v_{[0,n]} \in Q_1 \cap Q_2$, we then have also

$$\mathbb{P}(G^c | A \cap B, V_{[0,n]} = v_{[0,n]}) = \frac{\mathbb{P}(G^c \cap A \cap B | X_{[0,n]} \in \Xi, V_{[0,n]} = v_{[0,n]})}{\mathbb{P}(A \cap B | X_{[0,n]} \in \Xi, V_{[0,n]} = v_{[0,n]})}$$

$$\leq 2 e^{-(\nu - 2\alpha)n}$$

$$\leq \frac{1}{2},$$

which with (4.9) shows

$$(4.10) \qquad \mathbb{P}(G \cap A \cap B | X_{[0,n]} \in \Xi, V_{[0,n]} = v_{[0,n]}) \geq \tfrac{1}{4} e^{-(\tilde{I}(p_h \varepsilon) + \alpha) J^*}.$$

We claim that, for $\varepsilon$ sufficiently small,

$$(4.11) \qquad \tilde{I}(p_h \varepsilon) \leq 3 p_h I_V(\varepsilon).$$

To prove this, recall that $M_V(x) = \langle e^{xV_1} \rangle^V$ and let

$$f_1(x) = \psi(x) - p_h \varepsilon x,$$

$$f_2(x) = p_h \log M_V(x) - p_h \varepsilon x,$$

so that

$$\tilde{I}(p_h \varepsilon) = -\inf_x f_1(x),$$

$$p_h I_V(\varepsilon) = -\inf_x f_2(x).$$

The location $s_2(\varepsilon)$ of the infimum of $f_2$ satisfies

$$s_2(\varepsilon) \sim \frac{\varepsilon}{\mathrm{var}(V_1)}, \qquad f_2(s_2(\varepsilon)) \sim -\frac{p_h \varepsilon^2}{2\,\mathrm{var}(V_1)} \qquad \text{as } \varepsilon \to 0,$$

and for small $\varepsilon$,

$$t > \frac{3\varepsilon}{\operatorname{var}(V_1)} \qquad \text{implies } f_2(t) \geq \frac{\varepsilon^2 p_h}{2\operatorname{var}(V_1)}.$$

For fixed $v$, the function $\log(1 + p_h(e^{tv} - 1))$ is the log moment generating function of $v$ times a Bernoulli($p_h$) random variable, so it is convex in $t$; it follows that $\psi$ and $f_1$ are convex. Since $M_V$ is a moment generating function, $f_2$ is also convex. Also

$$f_1(0) = f_2(0) = 0,$$
$$f_1'(0) = f_2'(0) = -p_h\varepsilon.$$

It follows from all this that, to prove (4.11), it suffices to show that, for small $\varepsilon$,

$$\psi(x) \geq p_h \log M_V(x) - \frac{p_h\varepsilon^2}{2\operatorname{var}(V_1)} \qquad \text{for all } x \in \left[0, \frac{3\varepsilon}{\operatorname{var}(V_1)}\right],$$

and for this, in turn, it suffices to show

$$(4.12) \quad \psi(x) \geq \frac{19}{20} p_h \log M_V(x) - \frac{p_h\varepsilon^2}{4\operatorname{var}(V_1)} \qquad \text{for all } x \in \left[0, \frac{3\varepsilon}{\operatorname{var}(V_1)}\right],$$

since, for small $\varepsilon$,

$$\frac{1}{20} p_h \log M_V(x) \leq \frac{p_h\varepsilon^2}{4\operatorname{var}(V_1)}.$$

There exists $c_2$ such that, for all $x$,

$$|xV_1| \leq c_2 \qquad \text{implies } \log(1 + p_h(e^{xV_1} - 1)) \geq \tfrac{19}{20} p_h(e^{xV_1} - 1).$$

Define the event

$$C_x = \{|xV_1| \leq c_2\}.$$

Then using the concavity of log, for $x \in [0, 3\varepsilon/\operatorname{var}(V_1)]$,

$$\begin{aligned}
\psi(x) &= \langle \delta_{C_x} \log(1 + p_h(e^{xV_1} - 1)) \rangle^V \\
&\quad + \langle \delta_{C_x^c} \log(1 + p_h(e^{xV_1} - 1)) \rangle^V \\
&\geq \tfrac{19}{20} p_h \langle (e^{xV_1} - 1) \rangle^V - \tfrac{19}{20} p_h \langle \delta_{C_x^c}(e^{xV_1} - 1) \rangle^V \\
&\quad + \langle \delta_{C_x^c} \log(1 - p_h + p_h e^{xV_1}) \rangle^V \\
&\geq \tfrac{19}{20} p_h \log M_V(x) - \langle \delta_{C_x^c} e^{xV_1} \rangle^V + p_h x \langle \delta_{C_x^c} V_1 \rangle^V.
\end{aligned}$$

Since $V_1$ has exponential tails, there exists $c$ such that, provided $\varepsilon$ is small, for $x$ as above,

$$\langle \delta_{C_x^c} e^{xV_1} \rangle^V \leq e^{-c/x} \leq \frac{p_h \varepsilon^2}{8 \operatorname{var}(V_1)},$$

$$|p_h x \langle \delta_{C_x^c} V_1 \rangle^V| \leq e^{-c/x} \leq \frac{p_h \varepsilon^2}{8 \operatorname{var}(V_1)},$$

and (4.12) follows, proving (4.11).

Let $v_{[0,n]} \in Q_1 \cap Q_2 \cap Q_3$. From (4.11) and (4.10), provided $\alpha$ is small, since we chose $c_1 \leq 1/4$ in the definition of $J^*$, we have

$$(4.13) \qquad \mathbb{P}(G \cap A \cap B | X_{[0,n]} \in \Xi, V_{[0,n]} = v_{[0,n]}) \geq e^{-I_V(\varepsilon) p_g p_h \delta n}.$$

We also have

$$G \cap A \subset \{S_{\mathcal{R}} \geq (\varepsilon c_1 p_g p_h \delta - \alpha - \lambda) n\} \subset \left\{ S_{\mathcal{R}} \geq \frac{\varepsilon c_1 p_g p_h \delta}{2} n \right\},$$

which with (4.13) yields that, since $v_{[0,n]} \in Q_3$,

$$
\begin{aligned}
\mu_{[0,n]}^{\beta,u,\{v_i\}} & (L_n \leq \theta n) \\
& \leq \frac{\mu_{[0,n]}^{\beta,u,\{v_i\}}(L_n \leq \theta n)}{\mu_{[0,n]}^{\beta,u,\{v_i\}}(G \cap A \cap \{X_{[0,n]} \in \Xi\})} \\
& = \frac{\langle e^{\beta S_{\mathcal{R}}} e^{\beta u L_n} \delta_{\{L_n \leq \theta n\}} \rangle^X}{\langle e^{\beta S_{\mathcal{R}}} e^{\beta u L_n} \delta_{G \cap A \cap \{X_{[0,n]} \in \Xi\}}(\{v_i\}, \cdot) \rangle^X} \\
(4.14) \qquad & \leq \frac{h_\theta(v_{[0,n]}) \langle e^{\beta u L_n} \rangle^X}{e^{\beta \varepsilon c_1 p_g p_h \delta n/2} \langle e^{\beta u L_n} \delta_{G \cap A \cap \{X_{[0,n]} \in \Xi\}}(\{v_i\}, \cdot) \rangle^X} \\
& = \frac{h_\theta(v_{[0,n]})}{e^{\beta \varepsilon c_1 p_g p_h \delta n/2} \mathbb{P}(X_{[0,n]} \in \Xi) \mathbb{P}(G \cap A | X_{[0,n]} \in \Xi, V_{[0,n]} = v_{[0,n]})} \\
& \leq \mu_{[0,n]}^{\beta,u}(X_{[0,n]} \in \Xi)^{-1} \\
& \quad \times \exp\left( -\left( \frac{\beta \varepsilon c_1 p_g p_h \delta}{2} - I_V(\varepsilon) p_g p_h \delta - 2\beta \theta \right) n \right).
\end{aligned}
$$

Recall that $u = u_c^d - \eta$. For $\varphi(\delta)$ from Theorem 2.1, provided $\eta < \varphi(\delta)$, we have from that theorem that, for large $n$,

$$(4.15) \qquad \mu_{[0,n]}^{\beta,u}(L_n \geq \delta n) \geq e^{-\eta n} \mu_{[0,n]}^{\beta,u_c^d}(L_n \geq \delta n) \geq 2 e^{-2\varphi(\delta) n}.$$

By Lemma 2.3, we have

$$\mu_{[0,n]}^{\beta,u}(G_n \geq J^* | L_n \geq \delta n) \to 1 \qquad \text{as } n \to \infty.$$

With this we have from (4.14) and (4.15) that, for large $n$, for $v_{[0,n]} \in Q_1 \cap Q_2 \cap Q_3$,

$$
\begin{aligned}
(4.16) \quad & \mu_{[0,n]}^{\beta,u,\{v_i\}}(L_n \leq \theta n) \\
& \leq \exp\left(-\left(\frac{\beta \varepsilon c_1 p_g p_h \delta}{2} - I_V(\varepsilon) p_g p_h \delta - \beta a(\theta) - 2\varphi(\delta)\right) n\right).
\end{aligned}
$$

Since $I_V(\varepsilon) = O(\varepsilon^2)$ as $\varepsilon \to 0$, we can choose $\varepsilon$, then $\delta$, then $\theta, \eta$ small enough so that (4.16) implies that, for large $n$,

$$
\begin{aligned}
(4.17) \quad & \mu_{[0,n]}^{\beta,u,\{v_i\}}(L_n \leq \theta n) \leq \exp\left(-\frac{\beta \varepsilon c_1 p_g p_h \delta}{4} n\right) \\
& \hspace{4cm} \text{for all } \{v_i\} \in Q_1 \cap Q_2 \cap Q_3.
\end{aligned}
$$

It remains to show that $P_{[0,n]}^V(Q_2^c)$ is small. Let

$$
\xi_j = \delta_{\{j\text{th target is hit}\}}.
$$

Fix $\mathcal{J}$ with $|\mathcal{J}| = J^*$, fix $W^* = \{w_j, j \in \mathcal{J}\}$ and condition on the event $F = F(\mathcal{J}, W^*) = \{\mathfrak{G}^* = \mathcal{J}, \mathcal{W}^* = W^*\}$. Thus, $\mathcal{J}$ contains the indices of the good blocks, and $W^*$ contains the target locations for those blocks. The total reward received at hit targets is

$$
S_{\mathcal{W}} = \sum_{j \in \mathcal{J}} \xi_j V_{w_j}.
$$

The (conditional) log moment generating function of $S_{\mathcal{W}}$, given $F, V_{W^*}$, normalized by $J^*$, is

$$
\hat{\psi}(t | F, V_{W^*}) = \frac{1}{J^*} \sum_{j \in \mathcal{J}} \log(1 + p_h(e^{t V_{w_j}} - 1)),
$$

so

$$
\langle \hat{\psi}(t | F, V_{W^*}) \rangle^V = \psi(t).
$$

There exists $c_3$ such that

$$
|\psi''(t)| \leq c_3 \qquad \text{for all } |t| \leq \tfrac{1}{2}.
$$

Let $s_1 = s_1(\varepsilon)$ denote the location of the $f_1$ infimum, so that $\psi'(s) = p_h \varepsilon$, let

$$
a_V = \sup\{t \geq 0 : \langle e^{t|V_1|} \rangle^V < \infty\}
$$

and define the event

$$
\begin{aligned}
D_1 = D_1(W^*) = \Big\{ & v_{W^*} : |\hat{\psi}(t | F, v_{W^*}) - \psi(t)| \leq \frac{\alpha^2}{128 c_3} \\
& \text{for } t = s_1, s_1 - \frac{\alpha}{8 c_3} \text{ and } s_1 + \frac{\alpha}{8 c_3} \Big\}.
\end{aligned}
$$

Provided $\varepsilon, \alpha$ are small, we have $s_1 < \min(a_V, 1/4)$ and $s_1 + \alpha/8c_3 \le \frac{1}{2}$. When $v_{W^*} \in D_1$, we have, using convexity of $\hat{\psi}(\cdot|F, v_{W^*})$, that

$$
\begin{aligned}
(4.18) \qquad \hat{\psi}'(s_1|F, v_{W^*}) &\le \frac{\hat{\psi}(s_1 + (\alpha/(8c_3))|F, v_{W^*}) - \hat{\psi}(s_1|F, v_{W^*})}{\alpha/(8c_3)} \\
&\le \frac{\psi(s_1 + (\alpha/(8c_3))) - \psi(s_1) + \alpha^2/(64c_3)}{\alpha/(8c_3)} \\
&\le \psi'(s_1) + \frac{\alpha}{4}.
\end{aligned}
$$

Similarly,

$$
(4.19) \qquad \hat{\psi}'(s_1|F, v_{W^*}) \ge \psi'(s_1) - \frac{\alpha}{4}.
$$

Fix $v_{W^*} \in D_1$. Let

$$
u_j = \frac{p_h e^{s_1 v_{w_j}}}{1 - p_h + p_h e^{s_1 v_{w_j}}}
$$

and define probability measures

$$
m_j = p_h \delta_{v_{w_j}} + (1 - p_h)\delta_0, \qquad \tilde{m}_j = u_j \delta_{v_{w_j}} + (1 - u_j)\delta_0.
$$

Then $\tilde{m}_j$ is a tilted variant of $m_j$. If we change probabilities so that target $w_j$ is hit with probability $u_j$ (i.e., we consider $\{\xi_j v_{w_j} : j \in \mathcal{J}\}$ under the measure $\prod_{j \in \mathcal{J}} \tilde{m}_j$), then the mean of $(J^*)^{-1} S_{\mathcal{W}}$ (conditional on $\mathcal{J}, W^*$) becomes $\tilde{\mu}$ given by

$$
\tilde{\mu} = \hat{\psi}'(s_1|F, v_{W^*}) \in \left[ p_h \varepsilon - \frac{\alpha}{4}, p_h \varepsilon + \frac{\alpha}{4} \right],
$$

and the variance of $S_{\mathcal{W}}$ becomes

$$
\Phi(v_{W^*}) = \sum_{j \in \mathcal{J}} u_j(1 - u_j) v_{w_j}^2.
$$

Now since $s_1 < a_V$ and $u_j(1 - u_j) \le p_h e^{s_1 v_{w_j}}/(1 - p_h)$,

$$
\langle \Phi(V_{W^*}) \rangle^V \le \frac{p_h}{1 - p_h} \langle (V_1^2 e^{s_1 V_1}) \rangle^V J^* = c J^*.
$$

We define

$$
D_2 = D_2(W^*) = \left\{ v_{W^*} : \Phi(v_{W^*}) \le \frac{\alpha^2}{32} (J^*)^2 \right\}.
$$

Note that

$$
(4.20) \qquad |\tilde{\mu} J^* - \varepsilon c_1 p_g p_h \delta n| \le \frac{\alpha}{4} J^* + \varepsilon p_h \le \frac{\alpha}{2} J^*.
$$

For $v_{W^*} \in D_1 \cap D_2$, from Chebyshev's inequality, for large $n$,

$$(4.21) \qquad \left( \prod_{j \in \mathcal{J}} \tilde{m}_j \right) \left( |S_{\mathcal{W}} - \tilde{\mu} J^*| > \frac{\alpha}{4} J^* \right) \leq \frac{1}{2}.$$

Also, from the definition of $\tilde{m}_j$, for any bounded $A \subset \mathbb{R}$,

$$(4.22) \qquad \left( \prod_{j \in \mathcal{J}} \tilde{m}_j \right) (S_{\mathcal{W}} \in A) \leq \frac{e^{(\sup A) s_1 J^*}}{e^{\hat{\psi}(s_1 | F, v_{W^*}) J^*}} \left( \prod_{j \in \mathcal{J}} m_j \right) (S_{\mathcal{W}} \in A).$$

From (4.20)–(4.22),

$$\frac{1}{2} \leq \left( \prod_{j \in \mathcal{J}} \tilde{m}_j \right) \left( |S_{\mathcal{W}} - \tilde{\mu} J^*| \leq \frac{\alpha}{4} J^* \right)$$

$$\leq \frac{e^{s_1 (\tilde{\mu} + \alpha/4) J^*}}{e^{\hat{\psi}(s_1 | F, v_{W^*}) J^*}} \left( \prod_{j \in \mathcal{J}} m_j \right) (|S_{\mathcal{W}} - \varepsilon c_1 p_g p_h \delta n| \leq \alpha J^*),$$

and hence, for $v_{w^*} \in D_1 \cap D_2$, provided $\alpha$ is small,

$$\mathbb{P}(|S_{\mathcal{W}} - \varepsilon c_1 p_g p_h \delta n| \leq \alpha J^* | F, V_{W^*} = v_{W^*})$$

$$= \left( \prod_{j \in \mathcal{J}} m_j \right) (|S_{\mathcal{W}} - \varepsilon c_1 p_g p_h \delta n| \leq \alpha J^*)$$

$$(4.23) \qquad \geq \frac{1}{2} \exp\left( \hat{\psi}(s_1 | F, v_{W^*}) J^* - s_1 \left( \tilde{\mu} + \frac{\alpha}{4} \right) J^* \right)$$

$$\geq \frac{1}{2} \exp\left( \left( \psi(s_1) - s_1 p_h \varepsilon - \frac{\alpha^2}{128 c_3} - s_1 \alpha \right) J^* \right)$$

$$= \frac{1}{2} \exp\left( -\left( \tilde{I}(p_h \varepsilon) + \frac{\alpha^2}{128 c_3} + s_1 \alpha \right) J^* \right)$$

$$\geq \exp(-(\tilde{I}(p_h \varepsilon) + \alpha) J^*),$$

where the second inequality uses (4.20). In the other direction, since $s \leq 1/4$, provided $\alpha$ is small, using (4.20) again,

$$\mathbb{P}(|S_{\mathcal{W}} - \varepsilon c_1 p_g p_h \delta n| \leq \alpha J^* | F, V_{W^*} = v_{W^*})$$

$$\leq \left( \prod_{j \in \mathcal{J}} m_j \right) (\exp(s_1 S_{\mathcal{W}}) \geq \exp(s_1 (\varepsilon c_1 p_g p_h \delta n - \alpha J^*)))$$

$$(4.24) \qquad \leq \exp(\hat{\psi}(s_1 | F, v_{W^*}) J^* - s_1 (\varepsilon c_1 p_g p_h \delta n - \alpha J^*))$$

$$\leq \exp\left( \left( \psi(s_1) - s_1 p_h \varepsilon + \frac{\alpha^2}{128 c_3} + 2 s_1 \alpha \right) J^* \right)$$

$$\geq \exp(-(\tilde{I}(p_h \varepsilon) - \alpha) J^*).$$

Together (4.23) and (4.24) show that, for $|W^*| = J^*$,

$$(D_1(W^*) \cap D_2(W^*)) \times F(\mathcal{J}, W^*) \subset B.$$

Therefore,

$$
\begin{aligned}
(4.25) \qquad &\langle \mathbb{P}(B^c | X_{[0,n]} \in \Xi, V_{[0,n]}) \rangle^V \\
&= \mathbb{P}(B^c | X_{[0,n]} \in \Xi) \\
&\leq \sum_{\mathcal{J}} \sum_{W^*} P^V(D_1(W^*)^c \cup D_2(W^*)^c) P^X(F(\mathcal{J}, W^*) | \Xi) \\
&\leq \max_{W^*} P^V(D_1(W^*)^c \cup D_2(W^*)^c),
\end{aligned}
$$

where the sum and maximum are over $|\mathcal{J}| = J^*$ and $|W^*| = J^*$.

The random variables in the sum in the definition of $\hat{\psi}(t | F, V_{W^*})$ satisfy

$$\langle \exp(a \log(1 + p_h(e^{tV_1} - 1))) \rangle^V < \infty \qquad \text{for } |a| < 1, |t| < a_V,$$

and therefore there exists $\gamma_1, \rho_1 > 0$ such that

$$P^V(D_1(W^*)^c) \leq e^{-\gamma_1 J^*} \leq e^{-\rho_1 n}.$$

Let

$$U_j = \frac{p_h e^{s_1 V_{w_j}}}{1 - p_h + p_h e^{s_1 V_{w_j}}}, \qquad Y_j = U_j(1 - U_j) V_{w_j}^2.$$

Then $0 \leq Y_j \leq c_4$ for some $c_4 = c_4(p_h, s)$, so, for all $t$, we have $\langle \exp(tY_j) \rangle^V < \infty$. Therefore, there exist $\gamma_2, \rho_2 > 0$ such that

$$P^V(D_2(W^*)^c) \leq e^{-\gamma_2 J^*} \leq e^{-\rho_2 n},$$

and thus, by (4.25), for $\rho = \min(\rho_1, \rho_2)$,

$$\langle \mathbb{P}(B^c | X_{[0,n]} \in \Xi, V_{[0,n]}) \rangle^V \leq 2e^{-\rho n}.$$

Therefore, $P^V(Q_2^c) \leq 4e^{-\rho n}$.

With (4.7), (4.8) and (4.17), this shows that

$$\sum_{n=1}^{\infty} P^V\left(\mu_{[0,n]}^{\beta, u, \{V_i\}}(L_n \leq \theta n) > \exp\left(-\frac{\beta \varepsilon c_1 p_g p_h \delta}{4} n\right)\right) < \infty,$$

which, by the Borel–Cantelli lemma, shows that the polymer is pinned at $(\beta, u)$.

**5. Proofs of Theorem 1.6 and Corollary 1.7.**

PROOF OF THEOREM 1.6.   To prove (i), fix $M > 0$ and $n \geq 1$ to be specified, and fix a disorder realization $\{v_i\}$. We may assume the chain is aperiodic; there then exists $l_0$ such that $P^X(E_1 = l) > 0$ for all $l \geq l_0$. Let $l_1 = \lfloor \overline{G}_V(M)^{-1} \rfloor$ and assume $M$ is large enough so $l_1 \geq l_0$. Let $i_0 = 0$ and for $j \geq 1$, let $i_j = i_j(\{v_i\}) = \min\{i \in [i_{j-1} + l_0, i_{j-1} + l_1) : v_i \geq M\}$ if such an $i$ exists, and $i_j = i_{j-1} + l_1$ otherwise. Let

$$J_+ = J_+(\{v_i\}, n) = \{i_j : j \geq 1, i_j \leq n, v_{i_j} \geq M\},$$

$$J_- = J_-(\{v_i\}, n) = \{i_j : j \geq 1, i_j \leq n, v_{i_j} < M\},$$

$$J = J_+ \cup J_-, \qquad \theta = \theta(M) = \min\{P^X(E_1 = l) : l_0 \leq l \leq l_1\}.$$

We consider the set $\Upsilon_J$ of length-$n$ trajectories which return to 0 exactly at the times in $J \cap [1, n]$, and show that their Gibbs weight alone is enough to make the free energy positive. In fact, we have

$$
\begin{aligned}
Z_{[0,n]}^{\{v_i\}}(\beta, u) &\geq \sum_{\{x_i\} \in \Upsilon_J} \exp\left(\beta \sum_{i=1}^{n} (u + v_i)\delta_{\{x_i=0\}}\right) P_{[0,n]}^X(\{x_i, 0 \leq i \leq n\}) \\
&\geq \exp\left(\beta\left((u + M)|J_+| + \sum_{i \in J_-}(u + v_i)\right)\right) \\
&\quad\quad \times \prod_{k=1}^{|J|} P^X(E_1 = i_i - i_{k-1}) \cdot P^X(E_1 > i_{|J|+1} - i_{|J|}) \\
&\geq \exp\left(\beta\left((u + M)|J_+| + \sum_{i \in J_-}(u + v_i)\right)\right)\theta^{|J|+1}.
\end{aligned}
$$

(5.1)

For fixed $I$, conditionally on $\{J_-(\{V_i\}) = I\}$, the random variables $\{V_i, i \in I\}$ are i.i.d. with distribution $P^V(V_1 \in \cdot \,|\, V_1 < M)$. We may assume $M$ is large enough so that $E^V(V_1 | V_1 < M) > -1$. Then the events

$$\Psi_n = \left\{\{v_i\} : \sum_{i \in J_-(\{v_i\}, n)}(u + v_i) \geq (u - 1)|J_-(\{v_i\}, n)|\right\}$$

satisfy $P^V(\{V_i\} \in \Psi_n^c \text{ i.o.}) = 0$. Also, the random variables $i_j(\{V_i\}) - i_{j-1}(\{V_i\})$, $j \geq 1$, are i.i.d. with

$$
\begin{aligned}
P^V(i_j(\{V_i\}) \in J_-(\{V_i\})) &\leq P^V(i_j(\{V_i\}) - i_{j-1}(\{V_i\}) = l_1) \\
&= G_V(M)^{l_1 - l_0} \\
&< \tfrac{1}{2},
\end{aligned}
$$

the last inequality holding provided $l_1$ is large. Therefore, the events

$$\Phi_n = \{\{v_i\} : |J_-(\{v_i\}, n)| \le |J_+(\{v_i\}, n)|\}$$

satisfy $P^V(\{V_i\} \in \Phi_n^c \text{ i.o.}) = 0$. Note that $|J| + 1 \ge n/l_1$; it follows that, for $\{v_i\} \in \Psi_n \cap \Phi_n$, provided $M \ge 4|u| + 2$ and $\beta$ is large enough, we have by (5.1) that

$$(5.2) \quad Z_{[0,n]}^{\{v_i\}}(\beta, u) \ge \exp(\beta(2u + M - 1)|J|/2)\theta^{|J|+1} \ge \theta(\theta e^{\beta M/4})^{|J|} \ge e^{n/l_1},$$

so that $\beta f^q(\beta, u) \ge 1/l_1 > 0$. It then follows from Theorem 3.1 that there is pinning at $(\beta, u)$, so (i) is proved.

For (ii), we modify the above as follows. Fix $\beta > 0, n \ge 1$. For ease of exposition, we suppose for now that $G_V$ is continuous. Let

$$I_\beta = \left\{ l \ge k_0 : \min_{k_0 \le k \le l} P^X(E_1 = k) \ge e^{-\beta \overline{G}_V^{-1}(1/l)/8} \right\},$$

which is infinite by assumption. Let $l_0$ be the smallest element of $I_\beta$, let $l_1 \in I_\beta$ to be specified, and let $M = \overline{G}_V^{-1}(1/l_1)$. Since $G_V$ is continuous, we have $l_1 = \overline{G}_V(M)^{-1}$. Then (5.1), still holds, and we have

$$\theta \ge e^{-\beta \overline{G}_V^{-1}(1/l_1)/8} = e^{-\beta M/8}.$$

As in (5.2), we have for $\{v_i\} \in \Psi_n \cap \Phi_n$ that

$$Z_{[0,n]}^{\{v_i\}}(\beta, u) \ge \exp(\beta(2u + M - 1)|J|/2)\theta^{|J|+1} \ge \theta e^{\beta M|J|/8} \ge \theta e^{\beta Mn/16l_1},$$

and, provided $l_1$ is large, we have $P^V(\{V_i\} \in \Phi_n^c \cup \Psi_n^c \text{ i.o.}) = 0$, so that the free energy is again positive, and (ii) is proved.

In case $G_V$ is not continuous, there need not exist $M$ with $\overline{G}_V(M) = 1/l_1$, so we introduce auxiliary randomization. Let $\{U_i, i \ge 1\}$ be i.i.d. random variables, uniform in $[0, 1]$ and independent of $\{V_i\}$. We change $J_+$ by requiring that either $v_{i_j} > M$ or both $v_{i_j} = M$ and $U_j \le a$, with $a$ chosen so that

$$P^V(V_1 > M) + aP^V(V_1 = M) = 1/l_1.$$

Correspondingly, for $J_-$, we require that either $v_{i_j} < M$ or both $v_{i_j} = M$ and $U_j > a$. The rest of the proof requires only trivial modifications. $\square$

PROOF OF COROLLARY 1.7. From the remarks after the corollary statement, it is enough to prove (iii), and for this, it suffices to show that (1.5) and (1.6) imply that (1.4) holds for some $\{k_j\}$.

We may assume $p$ and $w$ are strictly monotone. Let $\varepsilon > 0$. We have by (1.6) that

$$\sum_k P^V((\zeta \circ p \circ w)^{-1}(\varepsilon V_1^+) \ge k) = \infty$$

and, therefore, by (1.5), there exists a sequence $\{m_j\}$ such that

$$
(5.3)
\begin{aligned}
\frac{1}{\lfloor w(m_j) \rfloor} &\le P^V((\zeta \circ p \circ w)^{-1}(\varepsilon V_1^+) \ge m_j) \\
&= \overline{G}_V\left(\frac{1}{\varepsilon} \log \frac{1}{p(w(m_j))}\right)
\end{aligned}
$$

for all $j$. Then letting $n_j = \lfloor w(m_j) \rfloor$, we have

$$
p(n_j) \ge p(w(m_j)) \ge e^{-\varepsilon \overline{G}_V^{-1}(1/n_j)} \qquad \text{for all } j.
$$

Since $\varepsilon$ is arbitrary, this shows there is a sequence $\{k_j\}$ satisfying (1.4). $\quad\square$

Department of Mathematics KAP 108
University of Southern California
Los Angeles, California 90089-2532
USA
E-mail: alexandr@usc.edu
URL: www-rcf.usc.edu/~alexandr

IMPA
Estrada Dona Castorina
Jardim Botânico CEP 22460-320
Rio de Janeiro
Brazil
E-mail: vladas@impa.br